\documentclass[11pt]{article}
\usepackage{amsfonts}

\usepackage{graphics}

\usepackage{indentfirst}
\usepackage{cite}
\usepackage{latexsym}
\usepackage{amsmath}
\usepackage{amssymb}
\usepackage[dvips]{epsfig}
\usepackage{amscd}
\def\on{\bar\rho}
\newtheorem{theorem}{Theorem}[section]
\newtheorem{remark}{Remark}[section]

\newtheorem{definition}{Definition}[section]
\newtheorem{lemma}[theorem]{Lemma}

\newtheorem{proposition}[theorem]{Proposition}

\newcommand{\vp}{\varphi}
\newcommand{\n}{\rho}

\newcommand{\ti}{\tilde}

\newcommand{\lm}{\lambda}

\renewcommand{\div}{ {\rm div }  }

\newcommand{\pa}{\partial}
\renewcommand{\r}{\mathbb{R}}
\newcommand{\bi}{\bibitem}
\renewcommand{\b}{B_{N_*}}

\newcommand{\ia}{\int_0^T}

\newcommand{\bt}{\begin{theorem}}
\newcommand{\bl}{\begin{lemma}}
\newcommand{\el}{\end{lemma}}
\newcommand{\et}{\end{theorem}}
\newcommand{\ga}{\gamma}

\newcommand{\al}{\alpha}
\newcommand{\de}{\delta}
\newcommand{\ve}{\varepsilon}
\newcommand{\la}{\label}

\newcommand{\ol}{\overline}

\newcommand{\bn}{\begin{eqnarray}}
\newcommand{\en}{\end{eqnarray}}
\newcommand{\bnn}{\begin{eqnarray*}}
\newcommand{\enn}{\end{eqnarray*}}

\newcommand{\bnnn}{\begin{eqnarray*}}
\newcommand{\ennn}{\end{eqnarray*}}

\newcommand{\ba}{\begin{aligned}}
\newcommand{\ea}{\end{aligned}}
\newcommand{\be}{\begin{equation}}
\newcommand{\ee}{\end{equation}}
\def\O{{\r^2 }}
\def\p{\partial}
\def\norm[#1]#2{\|#2\|_{#1}}
\def\lap{\triangle}

\def\o{\omega}

\newcommand{\no}{\nonumber\\}

\newcommand{\si}{\sigma}

\def\la{\label}

\def\na{\nabla}
\def\on{\bar\n}

\makeatletter      
\@addtoreset{equation}{section}
\makeatother

\makeatletter      
\@addtoreset{equation}{section}
\makeatother       

\title{ }

\date{}
\author{
 }

\title{Global Well-Posedness and Large Time Asymptotic Behavior of Classical Solutions to the   Compressible  Navier-Stokes Equations
  with  Vacuum }

\author{
  Jing L{\small I}$^{a,b},$ Zhouping  X{\small IN}$^b$\thanks{This research is supported in part by   Zheng Ge Ru Foundation, and Hong
Kong RGC Earmarked Research Grants CUHK-4041/11P, CUHK-4048/13P, a Focus Area Grant from
The Chinese University of Hong Kong, and a grant from Croucher Foundation.
  The research of \textsc{J. Li} was
partially supported  by the National Center for Mathematics and Interdisciplinary Sciences, CAS, and NNSFC  Grant No. 11371348.
 Email:   ajingli@gmail.com (J. Li),
 zpxin@ims.cuhk.edu.hk (Z. Xin).}
 \\
   {\normalsize  a.   Institute of Applied Mathematics, AMSS,}\\
   {\normalsize \&   Hua Loo-Keng Key Laboratory of Mathematics,} \\
  {\normalsize Chinese Academy of Sciences,} \\
   {\normalsize Beijing 100190,
 P. R. China  }\\
  {\normalsize   b.  The Institute of Mathematical Sciences,} \\
   {\normalsize  The Chinese University of Hong Kong, Hong
  Kong}}

\date{ }

\begin{document}
\maketitle

 \begin{abstract} This paper concerns the global well-posedness and large time asymptotic behavior of strong
 and classical solutions to the Cauchy problem of the
Navier-Stokes equations for viscous   compressible  barotropic flows  in two or three spatial dimensions with vacuum as far field density.  For strong and classical solutions, some a priori decay with rates (in large time) for both the pressure and the spatial gradient of  the velocity field are obtained  provided that the initial total energy is suitably {small.}   Moreover, by using these key decay rates  and some analysis on the expansion rates of the essential support of the density,   we  establish the global existence and uniqueness of classical
solutions (which may be of possibly large oscillations) in two  spatial dimensions, provided the smooth initial data are of small total energy.  In addition, the initial density can even have compact support. This, in particular, yields
the global regularity and uniqueness of the re-normalized weak solutions of Lions-Feireisl to the two-dimensional
compressible barotropic flows for all adiabatic number $\gamma>1$ provided that the initial total energy is small.
 \end{abstract}

Keywords: compressible Navier-Stokes equations;  global-wellposedness; large-time behavior; Cauchy problem; vacuum.

\section{Introduction}

We consider the  Navier-Stokes
equations
   \be \la{a1}  \begin{cases}\n_t+{\rm div} (\n u)=0,\\
 (\n u)_t+{\rm div}(\n u\otimes u)-\mu\Delta u-(\mu+\lambda)\na {\rm
div} u +\na P(\n) =0, \end{cases}\ee for viscous  compressible barotropic flows. Here, $t\ge 0$ is time,   $x\in \Omega\subset \r^N( N=2,3)$  is the spatial coordinate,  and $\n=\n(x,t),$
$u=(u^1,\cdots,u^N)(x,t),  $   and \be P(\n)=R\n^\ga\,\,( R>0, \ga>1)  \ee   are the fluid
density, velocity and pressure, respectively.  Without loss of generality, it is assumed that $R=1.$ The constant
viscosity coefficients $\mu$ and $\lambda$  satisfy the physical
restrictions: \be\la{h3} \mu>0,\quad 2\mu +  {N} \lambda\ge 0.
\ee
Let $\Omega=\r^N  $ and
 we consider  the Cauchy problem
for (\ref{a1}) with $(\n,u)$  vanishing at infinity (in some weak sense) with given initial data $\n_0$ and $u_0,$ as

\be \la{h2} \n(x,0)=\n_0(x), \quad \n u(x,0)=\n_0u_0(x),\quad x\in \Omega= \r^N.\ee

There are huge literatures on the large time existence and
behavior of solutions to (\ref{a1}). The one-dimensional problem
has been studied extensively, see
\cite{Kaz,Ser1,Ser2,Hof} and the references therein. For the
multi-dimensional case, the local existence and uniqueness of
classical solutions are known in \cite{Na,se1}  in the absence of
vacuum and recently, for strong solutions also, in \cite{K3,K1,
K2, S2,hlma} for the case that the initial density need not be positive
and may vanish in open sets. The global classical solutions were
first obtained by Matsumura-Nishida \cite{M1} for initial data
close to a non-vacuum equilibrium in some Sobolev space $H^s.$ In
particular, the theory requires that the solution has small
oscillations from a uniform non-vacuum state so that the density
is strictly away from vacuum. Later, Hoff \cite{H3,Hof2,Ho3}
studied the problem for discontinuous initial data. For the
existence of solutions for arbitrary data, the major breakthrough is due to Lions \cite{L1} (see also  Feireisl  \cite{F1,F2}), where the
global existence of weak solutions  when the exponent $\ga$ is suitably large are achieved. The
main restriction on initial data is that the initial total energy is
finite, so that the density vanishes at far fields, or even has
compact support. However, little is known on the structure of such
weak solutions, in particular, the regularity and the uniqueness of such weak solutions remain open.
This is a subtle issue, as  Xin  \cite{X1}  showed that in the case that the
initial density has compact support, any smooth solution in $C^1([0,T]:H^s(\r^d)) (s>[d/2]+2)$ to the
Cauchy problem of the full compressible Navier-Stokes
system without heat conduction blows up in finite time for any space
dimension $d\ge 1$, and the same holds for the isentropic case (\ref{a1}), at
least in one-dimension. The assumptions of \cite{X1} that
the initial density has compact support and that the smooth solution has finite energy are removed recently by Xin-Yan \cite{xy1} for a large class of initial data containing vacuum. However, this blow-up theory does not apply to
the isentropic flows in general, at least in the case of $\r^3$. Indeed, very recently, for the case that the initial density is allowed to vanish  and even has compact support, Huang-Li-Xin  \cite{hlx1} established the quite surprising
global existence and uniqueness of classical
solutions with constant state as far field which could be
either vacuum or non-vacuum to \eqref{a1}-\eqref{h2}  in three-dimensional space with smooth initial data which are of small total energy  but possibly large oscillations.
  Moreover,   it was also showed in \cite{hlx1} that for any $p>2,$
\be \la{poq1}\lim_{t\rightarrow \infty}\left(\|P(\n)-P(\ti\n)\|_{L^p(\r^3)}+\|\na u\|_{L^2(\r^3) }\right)=0,\ee where $\ti\n $
  is the constant far field density. This not only generalizes the classical results of Matsumura-Nishida \cite {M1},
but also yields the regularity and uniqueness of the weak solutions of Lions and Feireisl \cite{L1, F1, F2} with initial
data of small total energy.  Then a natural  question arises whether the theory of Huang-Li-Xin \cite{hlx1}
remains valid for the case of $\r^2$. This is interesting partially due to the following reasons: First, a positive
answer would yield immediately the regularity and uniqueness of weak solutions of Lions-Feireisl with small initial total energy whose existence has been proved for all $\ga>1$, see
 \cite{F1,F2}. Second, this question may be subtle due to the recent blow-up result in \cite{hlx3} where it is shown that non-trivial ͬtwo-dimensional spherically symmetric solution
in $C^1([0,T];H^s(\r^2) (s>2))$ with initial compactly supported density blows up in finite time. Technically, it is
not easy to modify the three-dimensional analysis of \cite{hlx1} to the two-dimensional case with initial density containing vacuum
since the analysis of \cite{hlx1} depends crucially on the a priori $L^6$-bound on the velocity. For two-dimensional
problems, only in the case that the far field  density is away from vacuum, the techniques of \cite{hlx1} can be modified directly since at this case, for any $p\in [2,\infty),$ the $L^p$-norm of a function $u$ can be bounded   by $\|\n^{1/2} u\|_{L^2}$  and $\|\na u\|_{L^2},$    and the similar results can be obtained (\cite{hlx3}). However, when the far field density is vacuum, it seems difficult to  bound the $L^p$-norm of $u$ by $\|\n^{1/2} u\|_{L^2}$  and $\|\na u\|_{L^2} $ for any $p\ge 1,$ so the global existence and large time behavior of strong or classical solutions to the Cauchy problem are much subtle and remain open. Therefore,  the main  aim  of this paper is to  study the  global  existence and large time behavior of   strong or classical solutions to   \eqref{a1}-\eqref{h2}  in some homogeneous  Sobolev spaces in two-dimensional space with vacuum as far field density,  and at the same time to investigate the decay rates of the pressure and the gradient of velocity in both two and three dimensional spaces provided the initial energy is suitably small, which turn  out to be one of the keys for the two-dimensional global well-posedness theory.

Before stating the main results, we first explain the notations and
conventions used throughout this paper. For $R>0$ and $\Omega=\r^N (N=2,3),$ set
$$B_R  \triangleq\left.\left\{x\in\Omega\right|
\,|x|<R \right\} , \quad \int fdx\triangleq\int_{\Omega }fdx.$$ Moreover, for $1\le r\le \infty, k\ge 1, $ and $\beta>0,$ the standard homogeneous and inhomogeneous Sobolev spaces are defined as follows:
   \bnn  \begin{cases}L^r=L^r(\Omega ),\quad D^{k,r}=D^{k,r}(\Omega)=\{v\in L^1_{\rm loc}(\Omega)| \na^k v\in L^r(\Omega)\}, \\ D^1=D^{1,2},\quad
W^{k,r}  = W^{k,r}(\Omega) , \quad H^k = W^{k,2} , \\ \dot H^\beta=\left\{f:\Omega
 \rightarrow \r\left|\|f\|^2_{\dot H^\beta}=
 \displaystyle{\int} |\xi|^{2\beta}|\hat f(\xi)|^2d\xi<\infty\right.
 \right\} ,\end{cases}\enn  where $\hat f$ is the Fourier transform
 of $f.$  Next, we also give the definition of strong solutions   as follows:
\begin{definition} If  all derivatives involved in \eqref{a1} for $(\rho,u)  $  are regular distributions, and   equations  \eqref{a1} hold   almost everywhere   in $\Omega\times (0,T),$ then $(\n,u)$  is called a  strong solution to  \eqref{a1}.
\end{definition}

For $\Omega=\r^N (N=2,3),$ the initial total energy  is defined as: \bnn  C_0 =
\int_{\Omega}\left(\frac{1}{2}\n_0|u_0|^2 + \frac{1}{\ga-1}P(\n_0) \right)dx. \enn  We consider first
the two-dimensional case, that is, $\Omega=\r^2.  $ Without loss of generality, assume that the initial density $\n_0$ satisfies
\be\la{oy3.7} \int_{\r^2} \n_0dx=1,\ee  which implies that there exists a positive constant $N_0$ such that  \be\la{oy3.8} \int_{B_{N_0}}  \n_0  dx\ge \frac12\int\n_0dx=\frac12.\ee

We can now state our first main result, Theorem \ref{th1}, concerning the global existence  of  strong solutions.
\begin{theorem}\la{th1} Let $\Omega=\r^2.$
 In addition to \eqref{oy3.7} and  \eqref{oy3.8},  suppose that  the initial data $(\n_0,u_0)$ satisfy  for any given numbers $M>0,$ $\on\ge 1,$   $a>1,q>2, $ and $\beta\in (0,1],$
\be \la{co1}  0\le   \n_0\le\bar{\n}, \quad  \bar x^a \rho_0\in   L^1 \cap H^1\cap W^{1,q},\quad u_0\in \dot
H^\beta \cap D^1 ,  \quad \n_0^{1/2}u_0  \in L^2,\ee
 and  \be\la{h7}
   \|u_0\|_{\dot H^\beta}+ \|\n_0\bar x^a\|_{L^1} \le M ,  \ee   where
  \be\la{2.07} \bar x\triangleq(e+|x|^2)^{1/2} \log^2 (e+|x|^2) .\ee Then there exists a positive constant $\ve$ depending
 on  $\mu ,  \lambda ,   \ga ,  a ,    \on, \beta, N_0,$ and $M$  such that if
 \be\la{i7}
     C_0\le\ve,
   \ee  the problem  \eqref{a1}-\eqref{h2} has a unique global strong solution $(\n,u)$ satisfying   for any $0<   T<\infty,$   \be\la{h8}
  0\le\n(x,t)\le 2\bar{\n},\quad  (x,t)\in \O\times[0,T],
  \ee \be\la{1.10}\begin{cases}
  \rho\in C([0,T];L^1 \cap H^1\cap W^{1,q} ),\\  \bar x^a\rho\in L^\infty( 0,T ;L^1\cap H^1\cap W^{1,q} ),\\ \sqrt{\n } u,\,\na u,\, \bar x^{-1}u,
   \,    \sqrt{t} \sqrt{\n}  u_t \in L^\infty(0,T; L^2 ) , \\ \na u\in  L^2(0,T;H^1)\cap  L^{(q+1)/q}(0,T; W^{1,q}), \\ \sqrt{t}\na u\in L^2(0,T; W^{1,q} )   ,  \\ \sqrt{\n} u_t, \, \sqrt{t}\na u_t ,\,  \sqrt{t} \bar x^{-1}u_t\in L^2(\O\times(0,T)) ,
   \end{cases}\ee and   \be \la{l1.2}\inf\limits_{0\le t\le T}\int_{B_{N_1(1+t)\log^\al(e+t) }}\n(x,t) dx\ge \frac14 ,\ee
for any $\al>1$  and  some positive constant $N_1$ depending only on $\al,N_0, $ and $M.$ Moreover, $(\n,u)$ has the following decay
rates,  that is,  for $t\ge 1,$   \be \la{lv1.2}\begin{cases}\|\na u(\cdot,t)\|_{L^p}  \le C(p)t^{-1+1/p}, \mbox{ for }p\in [2,\infty),\\  \|P(\cdot,t)\|_{L^r}\le C(r)t^{-1+1/r}, \mbox{ for }r\in (1,\infty),\\  \|\nabla \o(\cdot,t)\|_{L^2}+ \|\nabla F(\cdot,t)\|_{L^2}\le C t^{-1 }, \end{cases} \ee
where \be \la{hj1} \o\triangleq  \pa_1u^2-\pa_2u^1 ,\quad F\triangleq(2\mu+\lambda)\div u-P,\ee are respectively  the vorticity   and the effective viscous flux,  and $C(\alpha)$ depends on $\alpha$ besides  $  \mu ,  \lambda ,   \ga ,  a ,  \on, \beta, $ $ N_0,$ and $M.$
\end{theorem}

\begin{remark}  It should be noted here that the  decay rate estimates \eqref{lv1.2} combined with the estimate on
upper bound of the expansion rate of the essential support of the density  \eqref{l1.2}  play  a crucial role in deriving the  global existence of strong and classical solutions to the two-dimensional problem \eqref{a1}-\eqref{h2}. This is in contrast to the three-dimensional case (\cite{hlx1}) where the global existence of  classical solutions to  \eqref{a1}-\eqref{h2} was achieved without any bounds on the decay rates of the solutions partially due to the a priori $L^6$-bounds on the velocity field. As will be seen in the proof, the key observation is the decay with a rate
for the mean-square norm of the pressure in \eqref{lv1.2}.
\end{remark}

If the initial data $(\n_0,u_0)$ satisfy some additional regularity and compatibility conditions, the global strong  solutions obtained by Theorem \ref{th1}  become classical ones,   that is,
\begin{theorem}\la{t1} Let $\Omega=\r^2.$  In addition to the assumptions in Theorem \ref{th1},  assume further that  $(\n_0,u_0)$ satisfy
  \be\la{1.c1}
   \na^2 \n_0, \,\na^2 P(\n_0 )\in L^2\cap  L^q,  \quad  \bar x^{\de_0}\na^2   \rho_0  ,\,  \bar x^{\de_0}\na^2 P(  \rho_0 )     ,   \, \na^2 u_0 \in  L^2 ,
    \ee   for some     constant $\de_0\in (0,1), $  and the following  compatibility condition:
\be \la{co2}- \mu\lap u_0 - (\mu + \lm )\nabla \div u_0 +  \nabla P(\n_0)=\n_0^{1/2}g , \ee   with  some $g\in L^2 .$  Then,   in addition to \eqref{h8}-\eqref{lv1.2},    the  strong   solution  $(\rho,u)$ obtained by Theorem \ref{th1} satisfies  for
  any $0<  T<\infty,$  \be\la{1.a10}\begin{cases}
  \na^2\rho,   \,\,\na^2 P(\rho)\in C([0,T];L^2\cap L^q  ), \\ \bar x^{\de_0}\na^2  \rho   ,\,\,  \bar x^{\de_0} \na^2  P(  \rho),    \, \na^2 u   \in L^\infty( 0,T ;L^2 ) ,\\    \sqrt{\n}  u_t, \,\sqrt{t}   \na  u_t,\,\sqrt{t}  \bar x^{-1}  u_t,\, t\sqrt{\n}u_{tt}, \,t  \na^2 u_t\in L^\infty(0,T; L^2),\\ t\na^3 u\in  L^\infty(0,T; L^2\cap L^q), \,\\    \na u_t,\, \bar x^{-1}u_t,\,  t\na u_{tt},\, t\bar x^{-1}u_{tt}\in L^2(0,T;L^2), \\   t \na^2(\n u)\in L^\infty(0,T;L^{(q+2)/2}) .
   \end{cases}\ee    \end{theorem}

\begin{remark}  The solution obtained in Theorem \ref{t1} becomes
a classical one for positive time (\cite{hlma}). Although it has small energy, yet
whose oscillations could be arbitrarily large. In
particular, both
interior and far field vacuum are allowed.  There is
no requirement on the size of the set of vacuum states. Therefore, the initial density may have compact support. Moreover, by the strong-weak uniqueness theorem of Lions \cite{L1},
Theorem  \ref{th1} and Theorem \ref{t1} can be regarded as uniqueness and
 regularity
theory of Lions-Feireisl's weak solutions with
small initial energy, whose existence has been proved for all $\ga>1$ in \cite{L1,F2}.
\end{remark}

\begin{remark}
It is worth noting that the conclusions in
Theorem \ref{t1} and Theorem \ref{th1}  are somewhat surprising
since for the  isentropic compressible Navier-Stokes equations
 \eqref{a1}, any non-trivial two-dimensional
spherically symmetric solution $(\n,u)\in C^1([0,T],H^s) (s>2)$  with initial compact supported density blows up in finite time (\cite{hlx3}). Indeed, as in \cite{X1}, the key point  of   \cite{hlx3}  to prove the blowup phenomena  is based on the fact that the support of the density   will not grow in time in the space $C^([0,T];H^m).$    However, in
the current case, though the density $\n\in C([0,T];H^2),$ yet the velocity $u$ satisfies only  $\na u\in C((0,T];H^k).$ Note that the function    $u\in \{\na u\in H^k\}$ decays  much slower   for large values of the spatial variable $x$ than   $ u\in H^{k+1}.$ Therefore, it seems that it is the slow  decay of the velocity  field for large values of the spatial variable $x$   that leads to the global existence of smooth solutions. Unfortunately, such an argument cannot be valid
for the full compressible Navier-Stokes system since the blow-up results of Xin-Yan in \cite{xy1} work  for any classical
solutions with compactly supported initial density.
\end{remark}

For the three-dimensional case, that is,  $\Omega=\r^3,$   we have the following results concerning the decay properties of the global classical solutions whose existence is   essentially due to \cite{hlx1}.
\begin{theorem}\la{thv}  Let $\Omega=\r^3.$
  For given numbers $M>0,$ $\on\ge 1,$
 $\beta\in (1/2,1],$ and $q\in (3,6),$ suppose that
    the initial data $(\n_0,u_0)$ satisfy
\be \la{cvo1}  \rho_0 ,\, P(\n_0) \in H^2\cap W^{2,q},\quad P(\n_0),\,\n_0 | u_0|^2  \in L^1,\quad  u_0\in \dot
H^\beta,\quad \nabla u_0\in H^1 , \ee
 \be\la{hv7} 0\le \rho_0\le \bar{\rho},\quad
   \|u_0\|_{\dot H^\beta} \le M,   \ee
and the compatibility condition \be\la{cvo2} -\mu\triangle u_0
-(\mu+\lambda)\nabla \div u_0 + \nabla P(\rho_0) = \rho^{1/2}_0g, \ee
for some $ g\in L^2.$ Moreover, if $\ga>3/2,$ assume that \be \n_0\in L^1.\ee
 Then there exists a positive constant $\ve$ depending
 on $\mu,\lambda,   \ga, $  $\on, \beta,$ and $M$  such that if
 \be\la{iv7}
     C_0\le\ve,
   \ee  the Cauchy problem
  (\ref{a1})-(\ref{h2})
  has a unique global classical solution $(\rho,u)$ in
   $\r^3\times(0,\infty)$ satisfying for
  any $0<\tau<T<\infty,$
  \be\la{hv8}
  0\le\rho(x,t)\le 2\bar{\rho},\quad x\in \r^3,\, t\ge 0,
  \ee
   \be
   \la{hv9}\begin{cases}
   \rho  \in C([0,T];L^{3/2}\cap H^2\cap W^{2,q}),\\ P\in C([0,T];  L^1\cap H^2\cap W^{2,q}) , \quad  u\in L^\infty(0,T;L^6),\\
\na  u\in L^\infty(0,T;H^1)\cap L^2(0,T;H^2)\cap    L^\infty(\tau,T;H^2\cap W^{2,q}),\\
\na  u_t\in L^2(0,T;L^2) \cap
L^{\infty}(\tau,T;H^1)\cap H^1(\tau,T;L^2)  .\end{cases} \ee
 Moreover, for $r\in (1,\infty),$ there exist positive constants $C(r)$ and $C$ depending   on $\mu,\lambda,   \ga, $  $\on, \beta,$ and $M$ such that for $t\ge 1,$  \be \la{lvy8}\begin{cases}\|\na u(\cdot,t)\|_{L^p}  \le C t^{-1+1/p}, \mbox{ for }   p\in [2,6] , \\  \|P(\cdot,t)\|_{L^r}\le C(r)t^{-1+1/r},  \mbox{ for }   r\in (1,\infty),\\  \|\nabla (\nabla\times u)(\cdot,t)\|_{L^2}+ \|\nabla ((2\mu+\lambda)\div u-P )(\cdot,t)\|_{L^2}\le C t^{-1 }, \end{cases} \ee
 where if $\ga>3/2,$ $C(r)$ and $C$ both depend on $\|\n_0\|_{L^1(\r^3)}$ also.

\end{theorem}

\begin{remark} It should be pointed out that the large time asymptotic decay with rates of the global strong or classical solutions, \eqref{lv1.2} and \eqref{lvy8}, are completely new for the multi-dimensional compressible Navier-Stokes equations \eqref{a1} in the presence of vacuum. They show in particular that the $L^2$-norm of both the pressure and the gradient of the velocity decay in time with a rate $t^{-1/2}$, and the gradient of the vorticity  and the effective viscous flux decay faster than themselves. However,  whether the second derivatives of the velocity field decay  or not remains open.  This is an interesting problem and left for the future.

\end{remark}

We now make some comments on the analysis of this paper. Note that for initial
data in the class satisfying (\ref{co1}), (\ref{h7}), \eqref{1.c1}, and \eqref{co2} except $u_0\in
\dot H^\beta,$ the local existence and uniqueness of  classical
solutions to the Cauchy problem, (\ref{a1})-(\ref{h2}), have been
established recently in \cite{hlma}. Thus, to extend the classical
solution globally in time, one needs some global a priori estimates on
smooth solutions to (\ref{a1})-(\ref{h2}) in suitable higher norms.
  It turns out that as in the three-dimensional case \cite{hlx1}, the key issue here is to derive both the
time-independent upper bound for the density and the time-depending
higher norm estimates of the smooth solution $(\n, u)$, so some  basic ideas used in \cite{hlx1} will be adapted
here, yet new difficulties arises in the two-dimensional case. Indeed, the analysis in  \cite{hlx1} relies heavily on the basic fact that, for the three-dimensional case,  the $L^6$-norm of $v\in D^1(\r^3)$ can be bounded by the $L^2$-norm of the gradient of $v$  which fails for $v\in D^1(\r^2).$  In fact, for two-dimensional case,  some of the main new difficulties are due to the appearance of
vacuum at far field and the lack of integrability  of the velocity and its material derivatives  in the whole two-dimensional space. To overcome these difficulties, first, using the $L^1$-integrability of the density, we
observe that the $L^2$-norm in both space and time of the pressure   is time-independent (see \eqref{h27}).     This    is crucial to show that the $H^1$-norm of  the  effective viscous flux decays at the rate of $t^{-1/2}$ for large time (see \eqref{hg2}) which plays a key role in obtaining the decay property of the $L^\infty$-norm of the effective viscous flux. Then, after some careful estimates of the expansion rates  of the essential support of the density  (see \eqref{uq2}), we succeed in obtaining that,  for large time,   the $L^p$-norm of the gradient of the   effective viscous flux (see (\ref{hj1}) for the definition) can be  bounded by the product of $(1+t)^5$ and some function $g(t) $ whose temporal $L^2$-norm   is independent of time (see \eqref{z.2}).   Based on these key ingredients,  we
are able to obtain the desired estimates on $L^1(0,\min\{1,T\};\,L^\infty({\r}^2))$-norm and the time-independent ones on $L^4(\min\{1,T\},T;\,L^\infty({\r}^2))$-norm of the effective viscous flux (see \eqref{hg1}). Then, motivated by   \cite{lx}, we deduce from  these  estimates and Zlotnik's inequality (see
Lemma \ref{le1}) that the density admits a time-uniform upper
bound which is the key for global estimates of classical
solutions.  The
next main step is to bound the gradients of the density and the
velocity. Similar to  \cite{hlx1,hlx,hx2},  such bounds can be obtained by solving a logarithm
Gronwall inequality based on a Beale-Kato-Majda type inequality
(see Lemma \ref{le9}) and the a priori estimates we have just
derived, and  moreover, such a derivation yields simultaneously
also the bound for $L^1(0,T;L^\infty({\r}^2))$-norm of the
gradient of the velocity, see Lemma \ref{le4} and its proof.   Finally, with these a priori estimates on
the gradients of the density and the velocity  at hand, one can
estimate the higher order derivatives by using the same arguments
as in \cite{hx2,hlma} to obtain the desired results.

The rest of the paper is organized as follows: In Section \ref{se2}, we
collect some elementary facts and inequalities which will be needed
in later   analysis. Sections \ref{se3} and \ref{se5} are devoted to deriving the necessary
a priori estimates on classical solutions which are needed to extend
the local solution to all time. Then finally, the main results,
Theorems \ref{th1}-\ref{thv}, are proved in Section \ref{se4}.

\section{Preliminaries}\la{se2}

In this section, for $\Omega=\r^2,$ we will recall some  known facts and elementary
inequalities which will be used frequently later.

We begin with
the     local existence of strong and classical solutions whose proof   can be found in \cite{hlma}.

\begin{lemma}   \la{th0}  Let $\Omega=\r^2.$  Assume  that
 $(\n_0,u_0 )$ satisfies \eqref{co1}  except $u_0\in \dot H^\beta.$   Then there exist  a small time
$T >0$    and a unique strong solution $(\rho , u )$ to the
  problem   \eqref{a1}-\eqref{h2}  in
$\O\times(0,T )$ satisfying \eqref{1.10} and  \eqref{l1.2}. Moreover, if  $(\n_0,u_0)$ satisfies \eqref{1.c1} and    \eqref{co2} besides \eqref{co1}, $(\rho , u )$  satisfies \eqref{1.a10} also.
 \end{lemma}

Next, the following well-known Gagliardo-Nirenberg inequality (see \cite{nir})
  will be used later.

\begin{lemma}
[Gagliardo-Nirenberg]\la{l1} For  $p\in [2,\infty),q\in(1,\infty), $ and
$ r\in  (2,\infty),$ there exists some generic
 constant
$C>0$ which may depend  on $p,q, $ and $r$ such that for   $f\in H^1({\O }) $
and $g\in  L^q(\O )\cap D^{1,r}(\O), $    we have \be
\la{g1}\|f\|_{L^p(\O)}^p\le C \|f\|_{L^2(\O)}^{2}\|\na
f\|_{L^2(\O)}^{p-2} ,\ee  \be
\la{g2}\|g\|_{C\left(\ol{\O }\right)} \le C
\|g\|_{L^q(\O)}^{q(r-2)/(2r+q(r-2))}\|\na g\|_{L^r(\O)}^{2r/(2r+q(r-2))} .
\ee
\end{lemma}

The following weighted $L^p$ bounds for elements of the Hilbert space $  D^{1}(\O)  $ can be found in \cite[Theorem B.1]{L2}.
\begin{lemma} \la{1leo}
   For   $m\in [2,\infty)$ and $\theta\in (1+m/2,\infty),$ there exists a positive constant $C$ such that we have for all $v\in  D^{1,2}(\O),$ \be\la{3h} \left(\int_{\O} \frac{|v|^m}{e+|x|^2}(\log (e+|x|^2))^{-\theta}dx  \right)^{1/m}\le C\|v\|_{L^2(B_1)}+C\|\na v\|_{L^2(\O) }.\ee
\end{lemma}

The combination of Lemma \ref{1leo} with the Poincar\'{e} inequality yields

\begin{lemma}\la{lemma2.6} For $\bar x$   as in \eqref{2.07},
suppose that $\n  \in L^\infty(\O)$ is a   function such that
\be \la{2.12}  0\le \n\le M_1, \quad M_2\le \int_{\b}\n dx ,\quad \n \bar x^\alpha \in L^1(\r^2),\ee
for $ N_*\ge 1 $ and positive constants $   M_1,M_2, $  and   $\al.$  Then, for $r\in [2,\infty),$ there exists a positive constant $C$ depending only on $  M_1, M_2, \alpha,   $ and $ r$  such that
 \be\la{z.1}\left(\int_{\r^2}\n |v |^r dx\right)^{1/r}  \le C  N_*^3  (1+\|\n\bar x^\al\|_{L^1(\r^2)})  \left(  \|\n^{1/2} v\|_{L^2(\r^2)} + \|\na  v \|_{L^2(\r^2)}\right) ,\ee for each $v\in \left.\left\{v\in D^1 (\O)\right|\n^{1/2}v\in L^2(\r^2) \right\}.$

\end{lemma}

{\it Proof.} First, for $f\in L^1(\b),$ denote the average of $f$ over $\b$ by $$f_{\b}\triangleq \frac{1}{|\b|}\int_{\b}f(x)dx.$$  It then follows from \eqref{2.12} that
\be\la{rev2} \ba\left|\n_{\b}v_{\b}\right|=&  \left|\frac{1}{|\b|}\int_{\b}\left( \n_{\b}-\n\right)(v-v_{\b}) dx +\frac{1}{|\b|}\int_{\b}  \n v dx \right|\\ \le & 2M_1  N_*^{-1}\|v-v_{\b}\|_{L^2(\b)}+M_1^{1/2}N_*^{-1}\|\n^{1/2}v\|_{L^2(\b)} \\ \le & 8M_1  \|\na v \|_{L^2(\b)}+M_1^{1/2}N_*^{-1}\|\n^{1/2}v\|_{L^2(\b)},\ea\ee
where in the last inequality  one has used the following  Poincar\'{e} inequality (\cite[(7.45)]{la})
\be\la{lpa1}         \|v-v_{B_{N^*}}\|_{L^2(B_{N^*})}\le   4 N^*\|\na v\|_{L^2(B_{N^*}) }     . \ee

Then, it follows from \eqref{rev2}  and  \eqref{2.12} that
\bnn\ba \left|v_{\b}\right|\le C(M_1,M_2) N_*^2\|\na v \|_{L^2(\b) }+C(M_1,M_2) N_* \|\n^{1/2}v\|_{L^2(\b)},\ea\enn
which together with \eqref{lpa1} leads to
\be\la{3.15}  \ba    \int_{\b}|v|^2dx&\le 2  \int_{\b}|v-v_{\b}|^2dx+  2|\b||v_{\b}|^2 \\&\le   C(M_1,M_2) N_*^6\|\na v\|_{L^2(\b)}^2+C(M_1,M_2)N_*^{4}\|\n^{1/2}v\|_{L^2(\b)}^2.\ea\ee

Finally, it follows from Holder's inequality, \eqref{3h},    \eqref{3.15}, and   \eqref{2.12}  that for  $r\in [2,\infty)  $  and $\si=4/ (4+\alpha) \in (0,1),$
 \bnn\ba \int_{\r^2}\n |v |^r dx  &\le   \|(\n \bar x^\alpha)^\si\|_{L^{1/\si}(\r^2)} \| |v|^r\bar x^{-\alpha\si}\|_{L^{4/(\alpha\si)}(\r^2)} \|\n\|_{L^\infty(\r^2)}^{1-\si}\\ &\le C\left(1+\|\n \bar x^\al\|_{L^1(\r^2)}\right)   \left(  N_*^3\left(  \|\n^{1/2} v\|_{L^2(\r^2)}+ \|\na  v\|_{L^2(\r^2)}\right)\right)^r  ,\ea\enn
 which gives \eqref{z.1}. This completes the proof of Lemma \ref{lemma2.6}.

Next, for $ \nabla^{\perp}\triangleq (-\p_2,\p_1),$ denoting the
material derivative of $f $   by  $\dot f\triangleq
f_t+u\cdot\nabla f,$
we  state some elementary estimates which follow from (\ref{g1})
and the standard $L^p$-estimate  for the following elliptic system
derived from the momentum equations in (\ref{a1}): \be\la{h13}
\triangle F = \text{div}(\n\dot{u}),\quad\mu \triangle \o =
\nabla^\perp\cdot(\n\dot{u}) , \ee where $F$ and $\o$ are as in \eqref{hj1}.

\begin{lemma} \la{le3}
  Let $\Omega=\r^2$ and $(\n,u)$ be a smooth solution of
   (\ref{a1}).
    Then for   $p\ge 2$ there exists a   positive
   constant $C$ depending only on $p,\mu,$ and $\lambda$ such that
          \be\la{h19}
    \|{\nabla F}\|_{L^p(\O)} + \|{\nabla \o}\|_{L^p(\O)}
   \le C\norm[L^p(\O)]{\n\dot{u}},\ee \be
      \la{h20}\norm[L^p(\O)]{F} + \norm[L^p(\O)]{\o}
   \le C \norm[L^2(\O)]{\n\dot{u}}^{1-2/p }
   \left(\norm[L^2(\O)]{\nabla u}
   + \norm[L^2(\O)]{P }\right)^{2/p} ,
\ee \be \la{h18}
   \norm[L^p(\O)]{\nabla u} \le C \norm[L^2(\O)]{\n\dot{u}}^{1-2/p }
   \left(\norm[L^2(\O)]{\nabla u}
   + \norm[L^2(\O)]{P }\right)^{2/p}+
   C \norm[L^p(\O)]{P }.
  \ee
\end{lemma}
{\it Proof.} On the one hand, the standard $L^p$-estimate for the elliptic system
(\ref{h13}) yields (\ref{h19}) directly, which, together with
(\ref{g1}) and (\ref{hj1}), gives (\ref{h20}).
On the other hand, since $-\Delta u=-\na {\rm div}u -\na^\perp\o,$ we have \bn\la{kq1}\na u=-\na(-\Delta)^{-1}\na {\rm
div}u-\na(-\Delta)^{-1}\na^\perp \o.\en Thus applying the standard $L^p$-estimate to \eqref{kq1} shows  \bnn \ba \|\na u\|_{L^p(\O)}&\le C(p) (\|{\rm
div}u\|_{L^p(\O)}+\|\o\|_{L^p(\O)})\\ &\le C (p)  \norm[L^p(\O)]{F} +C(p)  \norm[L^p(\O)]{\o} +
   C(p)  \norm[L^p(\O)]{P },\ea  \enn which,
along with (\ref{h20}), gives (\ref{h18}). The proof of Lemma \ref{le3} is completed.

Next,   to get the
uniform (in time) upper bound of the density $\n,$ we need the following Zlotnik  inequality.
\begin{lemma}[\cite{zl1}]\la{le1}   Let the function $y$ satisfy
\bnn y'(t)= g(y)+b'(t) \mbox{  on  } [0,T] ,\quad y(0)=y^0, \enn
with $ g\in C(R)$ and $y,b\in W^{1,1}(0,T).$ If $g(\infty)=-\infty$
and \be\la{a100} b(t_2) -b(t_1) \le N_0 +N_1(t_2-t_1)\ee for all
$0\le t_1<t_2\le T$
  with some $N_0\ge 0$ and $N_1\ge 0,$ then
\bnn y(t)\le \max\left\{y^0,\overline{\zeta} \right\}+N_0<\infty
\mbox{ on
 } [0,T],
\enn
  where $\overline{\zeta} $ is a constant such
that \be\la{a101} g(\zeta)\le -N_1 \quad\mbox{ for }\quad \zeta\ge \overline{\zeta}.\ee
\end{lemma}

Finally,    the following Beale-Kato-Majda type inequality,
which was proved in \cite{bkm,kato} when $\div u\equiv 0,$   will be
used later to estimate $\|\nabla u\|_{L^\infty}$ and
$\|\nabla\n\|_{L^2\cap L^q} (q>2)$.
\begin{lemma}   \la{le9}  For $2<q<\infty,$ there is a
constant  $C(q)$ such that  the following estimate holds for all
$\na u\in L^2(\O)\cap D^{1,q}({\O }),$ \bnn \la{ww7}\ba \|\na
u\|_{L^\infty({\O })}&\le C\left(\|{\rm div}u\|_{L^\infty({\O })}+
\|\o\|_{L^\infty({\O })} \right)\log(e+\|\na^2
u\|_{L^q({\O })})\\&\quad+C\|\na u\|_{L^2(\O)} +C . \ea\enn
\end{lemma}

\section{\la{se3} A priori estimates(I): lower order estimates}

In this section, for $\Omega=\O,$ we will establish some necessary a priori bounds
for smooth solutions to the Cauchy problem (\ref{a1})-(\ref{h2}) to extend the local strong and classical solutions guaranteed by
Lemma \ref{th0}. Thus, let $T>0$ be a fixed time and $(\n,u)$ be
the smooth solution to (\ref{a1})-(\ref{h2})  on
${\r}^2 \times (0,T]$  with smooth initial
data $(\n_0,u_0)$ satisfying (\ref{co1}) and (\ref{h7}).

Set $\si(t)\triangleq\min\{1,t \}.$  Define
 \be\la{As1}
  A_1(T) \triangleq \sup_{   0\le t\le T  }\left(\sigma\|\nabla u\|_{L^2}^2\right) + \int_0^{T} \sigma\int
 \n|\dot{u} |^2 dxdt,
  \ee
  and   \be \la{As2}
  A_2(T)  \triangleq\sup_{  0\le t\le T   }\sigma^2\int\n|\dot{u}|^2dx + \int_0^{T}\int
  \sigma^2|\nabla\dot{u}|^2dxdt.
\ee

We have the following key a priori estimates on $(\n,u)$.
\begin{proposition}\la{pr1}  Under  the conditions of Theorem \ref{th1},
     there exists some  positive constant  $\ve$
    depending    on  $\mu ,  \lambda ,   \ga ,  a ,  \on, \beta,$ $N_0,$ and $M$  such that if
       $(\n,u)$  is a smooth solution of
       (\ref{a1})-(\ref{h2})  on $\O \times (0,T] $
        satisfying
 \be\la{z1}
 \sup\limits_{
 \O \times [0,T]}\n\le 2\bar{\n},\quad
     A_1(T) + A_2(T) \le 2C_0^{1/2},
   \ee
     the following estimates hold
        \be\la{z2}
 \sup\limits_{\O \times [0,T]}\n\le 7\bar{\n}/4, \quad
     A_1(T) + A_2(T) +\int_0^T\si \|P\|_{L^2}^2dt\le  C_0^{1/2},
  \ee
   provided $C_0\le \ve.$
\end{proposition}

The proof of Proposition \ref{pr1} will be postponed to the end of this section.

In the following, we will use the convention that $C$ denotes a
generic positive constant
 depending  on $\mu$, $\lambda$, $\gamma$, $a$,
$\bar{\n},$  $\beta,$ $N_0,$ and $M$, and  use $C(\al)$ to emphasize
that $C$ depends on $\al.$

We begin with the following   standard   energy estimate for
$(\n,u)$ and preliminary    $L^2$ bounds for $\nabla u$ and
$\n\dot{u}$.
\begin{lemma}\la{le2}
 Let $(\n,u)$ be a smooth solution of
 (\ref{a1})-(\ref{h2}) on $\O \times (0,T]. $
  Then there is a positive constant
  $C $ depending only  on $\mu,$  $\lambda,$  and $\gamma$  such that
  \be \la{a16} \sup_{0\le t\le T}\int\left(
\frac{1}{2}\n|u|^2+\frac{1}{\ga-1}P\right)dx + \mu \ia\int |\na
u|^2  dxdt\le C_0,\ee
  \be\la{h14}
  A_1(T) \le  C C_0 +C\sup\limits_{0\le t\le T}\|P\|_{L^2}^2 + C\int_0^{T}\sigma\int\left(|\nabla u|^3+P|\nabla u|^2\right)dx dt,
  \ee
 and
  \be\la{h15}
    A_2(T)
    \le   CA_1(T)  + C\int_0^{T} \sigma^2 \left(\|\nabla u\|_{L^4}^4+\|P\|_{L^4}^4\right) dt.
   \ee
\end{lemma}

{\it Proof.}
First,
the standard energy inequality reads:
\bnn
 \sup\limits_{0\le t\le T}\int\left(
\frac{1}{2}\n|u|^2+\frac{P}{\ga-1}\right)dx+\int_0^T \int\left( \mu |\na u|^2+(\mu+\lambda )(\div u)^2\right)dxdt \le  C_0,
\enn which together with \eqref{h3} shows \eqref{a16}.

Next,
multiplying $(\ref{a1})_2 $ by
$  \dot{u}  $ and then integrating the resulting equality over
${\O } $ lead  to \be\la{m0} \ba    \int  \n|\dot{u} |^2dx      &
= - \int  \dot{u}  \cdot\nabla Pdx + \mu \int \triangle
 u\cdot \dot{u}  dx + (\mu+\lambda)
 \int \nabla\text{div}u\cdot \dot{u}  dx . \ea \ee Since $P$ satisfies \be
\la{a95}P_t+u\cdot\nabla P+\ga P{\rm div}u=0 ,\ee
  integration by parts yields that \be\la{m1} \ba
  - \int    \dot{u} \cdot\nabla Pdx
= & \int (   (\text{div}u)_tP
-   (u\cdot\nabla u)\cdot\nabla P)dx  \\
= & \left(\int  \text{div}u P  dx \right)_t
  + \int  \left( (\ga-1)P
 (\text{div}u)^2+   P\p_iu_j\p_ju_i\right)dx \\
\le & \left(\int  \text{div}uPdx \right)_t
    +C \int P|\na u |^2  dx. \ea \ee
Integration by parts also implies that\be\la{m2} \ba
  \mu \int \triangle u\cdot \dot{u} dx
& = -\frac{\mu }{2}\left( \|\nabla u\|_{L^2}^2\right)_t
  -\mu   \int \p_iu_j\p_i(u_k\p_ku_j)dx \\
& \le  -\frac{\mu }{2}\left( \|\nabla u\|_{L^2}^2\right)_t  + C \int  |\nabla u|^3dx , \ea
\ee and that \be\la{m3}\ba (\mu+\lambda)
 \int \nabla\text{div}u\cdot \dot{u}  dx &= -
\frac{\lambda+\mu}{2}\left(
\|\text{div}u\|_{L^2}^2\right)_t - (\lambda+\mu)  \int\div
u\div(u\cdot\na u)dx \\& \le -\frac{\lambda+\mu}{2}\left(
\|\text{div}u\|_{L^2}^2\right)_t + C
\int  |\nabla u|^3dx  . \ea \ee
Putting \eqref{m1}-\eqref{m3} into  \eqref{m0}  leads to \be\la{n1}\ba
 B'(t) +\int  \n|\dot{u} |^2dx &\le
  C \int P|\na u |^2  dx
 + C  \|\nabla u\|_{L^3}^3  ,\ea\ee where \be \la{nv1}\ba
B(t)&\triangleq \frac{\mu  }{2}\|\nabla
u\|_{L^2}^2+\frac{ \lambda+\mu
}{2}\|\text{div}u\|_{L^2}^2-\int  \text{div}u P dx   \ea\ee
satisfies
\be \la{n2}\ba
\frac{\mu }{4}\|\nabla u\|_{L^2}^2 -C \|P\|_{L^2}^2 \le  B(t)
&\le C \|\nabla u\|_{L^2}^2+ C \|P\|_{L^2}^2 .  \ea\ee

Then, integrating \eqref{n1} multiplied by $\si$
over $(0,T)  $ and using \eqref{n2} and \eqref{a16} yield
(\ref{h14}) directly.

Finally, to prove \eqref{h15}, we will use the basic estimates of   $\dot u$ due to Hoff\cite{H3}. Operating $ \pa/\pa t+\div
(u\cdot) $ to $ (\ref{a1})_2^j ,$ one gets by some simple calculations that
\be\la{ax1}\ba &\n   (\dot u^j)_t+ \n u \cdot\na \dot u^j -\mu \lap \dot u^j-(\mu+\lm)\p_j(\div \dot u) \\&=\mu  \p_i(-\p_iu\cdot\na u^j+\div u\p_i u^j)-\mu  \div(\p_i u\p_i u^j)\\&\quad- (\mu+\lambda) \pa_j\left(\pa_i u\cdot\na u^i- (\div u)^2\right)-(\mu+\lambda) \div (\pa_ju\div u)\\&\quad+(\ga-1)\p_j(P\div u)+\text{div}(P\p_ju) .\ea\ee
Multiplying (\ref{ax1}) by $\dot u$  and integrating the resulting equation over $\O$ lead to\be\la{mm4} \ba  \left(
  \int\n|\dot{u}|^2dx \right)_t + { \mu} \int
 |\nabla\dot{u}|^2dx &
\le   C  \|\na
u\|_{L^4}^4 + C  \|P\|_{L^4}^4 ,\ea\ee
which  multiplied
  by $\si^2$ gives  (\ref{h15}) and  completes the proof of Lemma \ref{le2}.

\begin{remark}It is easy to check that the estimates \eqref{n1} and \eqref{mm4} also hold for $\Omega=\r^3.$\end{remark}

Next, we give a key observation that pressure decays in time.

\begin{lemma}\la{le5} Let $(\n,u)$ be a smooth solution  of
   (\ref{a1})-(\ref{h2})     on $\O \times (0,T] $ satisfying (\ref{z1}). Then there exists a positive constant $C(\on)$ depending only  on $\mu,$  $\lambda,$  $\gamma,$  and $\on$
 such that
  \be\la{h27}
  A_1(T)+A_2(T)+\int_0^T\si\|P\|_{L^2}^2dt\le C(\on) C_0 .
  \ee
   \end{lemma}

{\it Proof.} First, it follows from (\ref{h18}),   \eqref{a16},  and   \eqref{z1} that
  \be\la{h99} \ba
  &  \int_0^{T}\sigma^2  \left(\|\na u\|_{L^4}^4 +\|P\|_{L^4}^4 \right) dt\\
& \le  C \int_0^{T}\sigma  \|\n  \dot u \|_{L^2}^2\left(\si\|\na u\|_{L^2}^2+\si\|P\|_{L^2}^2\right)dt+C  \int_0^{T}\sigma^2  \|P\|_{L^4}^4  dt\\
& \le  C(\on)\left(A_1 (T)+C_0 \right)\int_0^{T}\sigma  \|\n^{1/2}  \dot u \|_{L^2}^2dt+C (\on) \int_0^{T}\sigma^2  \|P\|_{L^2}^2   dt  . \ea \ee
To estimate the last term on the right-hand side of \eqref{h99}, noticing that $(\ref{a1})_2$ gives
\be\la{u} P=(-\Delta )^{-1}\div(\n  \dot u )+(2\mu+\lambda)\div u,\ee
we obtain from  H\"older's  and Sobolev's inequalities   that
\bnn \ba \int P^2dx \le&C \|(-\Delta)^{-1} \div(\n  \dot u )\|_{L^{4\ga}} \|P\|_{L^{4\ga/(4\ga-1)}} +C\|\na u\|_{L^2}\|P\|_{L^2} \\ \le& C\|  \n  \dot u \|_{L^{4\ga/(2\ga+1)}}\|\n\|_{L^1}^{1/2}\|\n\|_{L^{2\ga}}^{\ga-1/2}+C\|\na u\|_{L^2} \|P\|_{L^2} \\ \le& C \|\n^{1/2}\|_{L^{4\ga}} \| \n^{1/2}  \dot u  \|_{L^2}\|\n\|_{L^1}^{1/2}\|\n\|_{L^{2\ga}}^{\ga-1/2} +C\|\na u\|_{L^2} \|P\|_{L^2} \\ \le& C\|P\|_{L^2}  \| \n^{1/2}  \dot u  \|_{L^2} +C\|\na u\|_{L^2} \|P\|_{L^2} ,\ea \enn where in the last inequality, one has used  \be \la{mr3}\int\n dx =\int \n_0dx =1,\ee due to   the mass conservation equation  $\eqref{a1}_1.$
Thus, we arrive at
\be \la{new1}\ba \|P\|_{L^2}  \le   C  \| \n^{1/2}  \dot u \|_{L^2} + { C}\|\na u\|_{L^2}  ,\ea \ee
which, along with  \eqref{h14},   \eqref{h15}, \eqref{h99}, \eqref{a16}, and \eqref{z1} gives   \be\la{h28}\ba
A_1(T)+A_2(T)\le& C(\on)C_0 +C(\on)\int_0^{T} \sigma\|\nabla u\|_{L^3}^3dt. \ea\ee

Then, on the one hand,  one deduces from (\ref{h18}),   (\ref{a16}),  and (\ref{z1})  that \be\la{h34} \ba
  \int_0^{\si(T)} \sigma\|\nabla u\|_{L^3}^3 dt
&  \le C \int_0^{\si(T)} \sigma\|\n^{1/2} \dot u \|_{L^2} \left(\|\na u\|_{L^2}^2+\|P\|_{L^2}^2 \right)  dt +C(\on)C_0\\
& \le C A_2^{1/2}(\si(T))\int_0^{\si(T)}  \left(\|\na u\|_{L^2}^2+\|P\|_{L^2}^2 \right)  dt +C(\on)C_0 \\
& \le C(\on) C_0  . \ea \ee
  On the other hand, H\"older's inequality, \eqref{h99},  \eqref{z1}, and \eqref{new1} imply \be\la{h33} \ba
\int_{\si(T)}^{T} \sigma\|\nabla u\|_{L^3}^3 dt &\le
\de\int_{\si(T)}^{T}  \|\nabla u\|_{L^4}^4 dt + C(\de) \int_{\si(T)}^{T}\|\nabla u\|_{L^2}^2 dt \\& \le  \de C(\on) A_1(T)+  C(\de)C(\on) C_0.
\ea \ee

 Finally, putting   \eqref{h34}  and \eqref{h33} into \eqref{h28} and choosing $\de$ suitably small lead to \bnn A_1(T)+A_2(T)\le C(\on)C_0,\enn  which together with \eqref{new1} and \eqref{a16} gives \eqref{h27} and completes  the proof of Lemma \ref{le5}.

Next, we derive the rates of decay for $\na u$ and $P$, which are  essential to obtain the uniform  (in time) upper bound of the density for large time.
\begin{lemma}\la{ly1}
    For $p\in [2,\infty),$  there exists a positive constant $C(p,\on)$
  depending only  on $p, \mu,$  $\lambda,$  $\gamma,$  and $\on$ such that,  if  $(\n,u)$ is a smooth solution  of
   (\ref{a1})-(\ref{h2})     on $\O \times (0,T] $ satisfying (\ref{z1}), then
\be \la{ly8} \ba & \sup\limits_{
\si(T)\le t\le T}\left(t^{p-1}(
 \|\na u\|_{L^p}^p  +  \| P\|_{L^p}^p  )+t^2\|\n^{1/2}\dot u\|_{L^2}^2 \right) \le C(p,\on)C_0.\ea\ee
  \end{lemma}

{\it Proof.} First,  for $p\ge 2,$ multiplying (\ref{a95}) by $p P^{p-1}$ and integrating
the resulting equality over ${\O } ,$ one gets after using ${\rm
div}u=\frac{1}{2\mu+\lambda}(F+P)$ that \be\la{vv1}\ba
\left(\| P\|_{L^p}^p \right)_t+ \frac{p\ga-1}{2\mu+\lambda}\|P\|_{L^{p+1}}^{p+1} &
=- \frac{p\ga-1}{2\mu+\lambda}\int
P^pFdx \\&
\le \frac{p\ga-1}{2(2\mu+\lambda)}\|P\|_{L^{p+1}}^{p+1}+C(p)  \|F\|_{L^{p+1}}^{p+1} ,\ea\ee
which together with \eqref{h20} gives
\be\la{a96}\ba
\frac{2(2\mu+\lambda)}{p\ga-1}\left(\| P\|_{L^p}^p \right)_t+ \|P\|_{L^{p+1}}^{p+1}
    &\le C(p) \|F\|_{L^{p+1}}^{p+1} \\ &\le C(p) \left(\|\na u\|_{L^2}^2+\|P\|_{L^2}^2\right) \|\n \dot u\|_{L^2}^{p-1} .\ea\ee
In particular, choosing $p=2$  in \eqref{a96} shows
\be\ba
\la{a9z6}\left(\| P\|_{L^2}^2 \right)_t+ \frac{2\ga-1}{2(2\mu+\lambda)}\|P\|_{L^3}^3  \le
 \de \|\n^{1/2}\dot u\|_{L^2}^2+C(\de)\left(\|\na u\|_{L^2}^4+\|P\|_{L^2}^4\right) .\ea\ee

Next, it follows from \eqref{n1} and \eqref{h18} that
\be   \la{ly11}\ba
  B' (t)  + \int  \n|\dot{u} |^2dx & \le
   C  \| P\|_{L^3}^3
 + C  \|\nabla u\|_{L^3}^3   \\& \le
  C_1  \|P\|_{L^3}^3+ C \|\n\dot u\|_{L^2}\left(\|\na u\|_{L^2}^2+\|P\|_{L^2}^2\right) \\& \le
  C_1 \|P\|_{L^3}^3+ \de\|\n^{1/2}\dot u\|_{L^2}^2+C(\on,\de)\left(\|\na u\|_{L^2}^4+\|P\|_{L^2}^4\right).\ea\ee
 Choosing $ C_2\ge 2+  2(2\mu+\lm)(C_1+1)/(2\ga-1)$
suitably large such that \be\la{ly12} \frac{\mu}{4}\|\na u\|_{L^2}^2+\|P\|_{L^2}^2\le B(t)+C_2\|P\|_{L^2}^2\le C\|\na u\|_{L^2}^2+C\|P\|_{L^2}^2,\ee  adding \eqref{a9z6} multiplied by $C_2$ to \eqref{ly11},  and choosing $\de$ suitably small lead to
\be \la{h29}\ba &
 2\left( B (t)+  C_2 \| P\|_{L^2}^2\right)' +    \int  \left(\n|\dot{u} |^2 + P^3  \right)dx  \le
   C \| P\|_{L^2}^4 + C  \|\na u\|_{L^2}^4  ,\ea\ee which  multiplied  by $t,$  together with Gronwall's inequality,   \eqref{ly12},    \eqref{h27},    \eqref{a16},  and  \eqref{z1}    yields
 \be \la{h31} \ba & \sup\limits_{
\si(T)\le t\le T}t\left( \| \na u\|_{L^2}^2+ \| P\|_{L^2}^2
 \right)  + \int_{\si(T)}^Tt \int  \left(\n|\dot{u} |^2 + P^3  \right)dxdt\le C(\on)C_0. \ea\ee

Next,
multiplying \eqref{mm4} by $t^2 $ together with \eqref{h18} gives \be\la{zo1} \ba & \left(t^2
  \int\n|\dot{u}|^2dx \right)_t + { \mu}t^2 \int
 |\nabla\dot{u}|^2dx \\&
\le   2t   \int\n|\dot{u}|^2dx+ C (\on)t^2\|\n\dot u\|_{L^2}^2\left(
 \|\na u\|_{L^2}^2+ \|P \|_{L^2}^2\right) + \ti C (\on)t^2\|P\|_{L^4}^4.\ea\ee
     Choosing $p=3$ in \eqref{a96} and adding \eqref{a96}  multiplied by $(\ti C+1)t^2$ to \eqref{zo1} lead  to
\bnn\la{ly7} \ba & \left(t^2
  \int\n|\dot{u}|^2dx +  \frac{2(2\mu+\lambda)(\ti C +1)}{3\ga-1}t^2\| P\|_{L^3}^3\right)_t + { \mu}t^2 \|\nabla\dot{u}\|^2_{L^2}+t^2\|P\|_{L^4}^4   \\&
\le   C   t \int\left(\n|\dot{u}|^2+   P^3\right)dx+ C (\on)t^2\|\n^{1/2} \dot u\|_{L^2}^2 \left(\|\na u
 \|_{L^2}^2+  \|P
 \|_{L^2}^2\right),\ea\enn
 which combined with  Gronwall's inequality,     \eqref{h31},  and   \eqref{z1} yields
\be \la{ly9} \ba   \sup\limits_{\si(T)\le t\le T}t^2\int  \left(\n|\dot{u} |^2 + P^3  \right)dx +  \int_{\si(T)}^Tt^2 \left(\|\nabla\dot{u}\|^2_{L^2}+ \|P\|_{L^4}^4\right)dt \le C(\on)C_0. \ea\ee

Finally, we claim that for   $m=1,2,\cdots,
$ \be \la{ly10}\sup\limits_{\si(T)\le t\le   T}t^{m}\|P\|_{L^{m+1}}^{m+1}+\int_{\si(T)}^Tt^m\|P\|_{L^{m+2}}^{m+2}dt\le C(m,\on)C_0,\ee which together with \eqref{h18}, \eqref{h31}, and \eqref{ly9} gives \eqref{ly8}.  We shall prove \eqref{ly10} by induction. In fact, \eqref{h31} shows that \eqref{ly10} holds for $m=1.$ Assume that \eqref{ly10} holds for $m=n,$ that is, \be \la{ly16}\sup\limits_{\si(T)\le t\le T}t^{n}\|P\|_{L^{n+1}}^{n+1}+\int_{\si(T)}^Tt^n\|P\|_{L^{n+2}}^{n+2}dt\le C(n,\on)C_0.\ee
Multiplying \eqref{a96} where  $p=n+2$  by $t^{n+1} ,$ one obtains after  using \eqref{ly9}
\be\la{zo8}\ba &
  \frac{2(2\mu+\lambda)}{(n+2)\ga-1}\left(t^{n+1}\| P\|_{L^{n+2}}^{n+2} \right)_t+ t^{n+1}\|P\|_{L^{n+3}}^{n+3} \\ &
\le C (n,\on) t^{n}\| P\|_{L^{ n+2} }^{n+2} +C(n,\on) C_0 \left(\|\na u\|_{L^2}^2+\|P\|_{L^2}^2\right).\ea\ee Integrating \eqref{zo8} over $[\si(T),T]$ together with \eqref{ly16} and  \eqref{h27} shows that \eqref{ly10} holds for $m=n+1.$   By induction, we obtain \eqref{ly10} and finish the proof of Lemma \ref{ly1}.

 Next, the following Lemma \ref{lemma2.7} combined with Lemma \ref{lemma2.6}  will be useful to estimate the $L^p$-norm of $\n \dot u$   and obtain
     the uniform  (in time)
upper bound of the density for large time.

\begin{lemma}\la{lemma2.7} Let $(\n,u)$ be a smooth solution of (\ref{a1})-(\ref{h2}) on $\O \times (0,T] $ satisfying the assumptions in Theorem \ref{th1} and (\ref{z1}). Then for any $\al>0,$ there
exists a    positive constant  $N_1$   depending only on   $\al,$     $N_0,$ and $M$ such
that for all $t\in (0,T],$
    \be\la{uq2}   \int_{B_{N_1(1+t)\log^\al(1+t)}}\n(x,t) dx \ge  \frac14. \ee

\end{lemma}

{\it Proof.}  First, multiplying $\eqref{a1}_1$ by $(1+|x|^2)^{1/2} $ and integrating the resulting equality over $\O,$ we obtain after integration by parts and using both \eqref{a16} and  \eqref{mr3} that
\bnn \ba  \frac{d}{dt}\int\n (1+|x|^2)^{1/2}  dx &\le C\int \n |u|  dx\\ &\le C\left(\int\n   dx\right)^{1/2}\left(\int\n |u|^2 dx\right)^{1/2}\\ &\le C.\ea\enn
  This  gives
\be \la{o3.7} \sup_{0\le s\le t}\int\n  (1+|x|^2)^{1/2}  dx\le  C(M)(1+t) .\ee

 Next,
for $ \varphi(y)\in C^\infty_0(\r^2) $   such that \bnn 0\le\varphi(y)\le 1,\quad \varphi(y) =\begin{cases} 1& \mbox{ if }\,|y|\le 1,\\ 0& \mbox{ if }\,|y|\ge 2 ,\end{cases} \quad |\na \varphi|\le 2,\enn multiplying $\eqref{a1}_1$ by $\vp(y)$ with
  $y= \ti \de x(1+t)^{-1}\log^{-\al}(e+t)$   for   small $\ti \de>0$ which will be determined later, we obtain  \bnn\ba \frac{d}{dt}\int \n \vp(y)dx&= \int\n \na_y\vp\cdot y_tdx+\frac{\ti\de}{(1+t)\log^\al(e+t)  } \int \n u \cdot\na_y \vp  dx\\&\ge  - \frac{C\ti\de}{(1+t)^2\log^\al(e+t)  } \int\n| x |  dx- \frac{C\ti\de}{(1+t)\log^\al(e+t)  }  \\ &\ge -  \frac{C(M)\ti\de}{(1+t)\log^\al(e+t)  } ,\ea\enn where in the last inequality we have used \eqref{o3.7}.
Since $\al>1,$ this yields
\be\la{oi1} \ba \int\n\vp\left(y\right)dx&\ge \int \n_0(x)\vp(x\ti\de)dx-C(\al,M)\ti\de\ge \frac14 ,\ea\ee where we choose $\ti \de=(N_0+4C(\al,M))^{-1}. $

Finally, it follows from \eqref{oi1} that  for   $N_1\triangleq 2\ti \de^{-1}=2(N_0+4C(\al,M)),$ \bnn\ba \int_{B_{N_1(1+t)\log^\al(e+t)}}\n dx&\ge \int\n\vp\left( \ti \de x(1+t)^{-1}\log^{-\al}(e+t)\right)dx\ge \frac14, \ea\enn which  finishes the proof of Lemma \ref{lemma2.7}.

Next,    to  obtain      the
upper bound of the density for small time, we still need the following lemma.

\begin{lemma}\la{zc1} Let $(\n,u)$ be a smooth solution of (\ref{a1})-(\ref{h2})   on $\O \times (0,T] $ satisfying (\ref{z1}) and the assumptions in Theorem \ref{th1}. Then there
exists a positive constant  $K $   depending only on $\mu ,  \lambda ,   \ga ,  a ,  \on, \beta,$ $N_0,$ and $M$  such
that
   \be\la{uv1}  \sup_{0\le t\le \si(T)}t^{1-\beta}\|\na
u\|_{L^2}^2+\int_0^{\si(T)}t^{1-\beta}\int\n|\dot u|^2dxdt\le
K(\on,M). \ee

\end{lemma}

{\it Proof.} First, set $$\nu\triangleq \min\left\{\frac{\mu^{1/2}}{ 2(1+2\mu+\lm)^{1/2} },\,\frac{ \beta}{1-\beta}\right\}\in (0,1/2]  .$$ If $\beta\in (0,1),$   Sobolev's inequality implies
\be \la{zz.1}\ba\int \n_0|u_0|^{2+\nu}dx &\le  \int \n_0|u_0|^{2 }dx +  \int \n_0 |u_0|^{2/(1-\beta)}dx\\ &\le C(\on)+C(\on)\|u_0\|_{\dot H^\beta}^{2/(1-\beta)}\le C(\on,M).\ea\ee For the case that $\beta=1,$ one obtains from \eqref{z.1} that
\be  \la{zz.2}\ba\int \n_0|u_0|^{2+\nu}dx &\le C(\on)\left( \int \n_0|u_0|^{2 }dx +  \int   |\na u_0|^2dx\right)^{(2+\nu)/2} \le C(\on,M).\ea\ee

Then, multiplying $(\ref{a1})_2$ by $(2+\nu)|u|^\nu u$  and   integrating the resulting equation over $ \O$ lead  to
\bnn\ba & \frac{d}{dt}\int \n |u|^{2+\nu}dx+ (2+\nu) \int|u|^\nu \left(\mu |\na u|^2+(\mu+\lm) (\div u)^2\right) dx \\& \le  (2+\nu)\nu  \int (\mu+\lm)|\div u||u|^\nu |\na u|dx +C  \int \n^\ga |u|^\nu |\na u|dx \\& \le \frac{2+\nu}{2}\int (\mu+\lm)(\div u)^2|u|^\nu  dx +\frac{(2+\nu)\mu}{4} \int  |u|^\nu |\na u|^2 dx \\&\quad + C \int \n |u|^{2+\nu} dx+C\int \n^{(2+\nu)\ga-\nu/2}dx ,\ea\enn
which  together with Gronwall's inequality, \eqref{zz.1}, and \eqref{zz.2}   thus gives
   \be\la{jan2}
  \sup_{0\le t\le  \si(T) }\int \n |u|^{2+\nu}dx\le C(\on,M ) .
  \ee

Next, as in \cite{hof2002},  for  the linear differential operator $L$
defined by \be\ba(Lw)^j&\triangleq \n w^j_{t}+ \n u\cdot\na w^j
-(\mu\Delta w^j+(\mu+\lambda) \pa_j\div w )\no&=\n\dot w^j
-(\mu\Delta w^j+(\mu+\lambda) \pa_j\div w ),\quad j=1,2  , \ea\ee let $w_1$ and $w_2$   be
the solution to: \be \la{sas2}Lw_1 =0 ,\quad w_1(x,0)=w_{10}(x), \ee
and\be\la{sas3} Lw_2 = -\nabla P(\n),\quad w_2(x,0)=0,\ee
respectively.
A straightforward energy estimate of \eqref{sas2}  shows  that: \be\la{sas4}
\sup_{0\le t\le\si(T)}\int \n|w_1|^2dx+\int^{\si(T)}_{0}
\int|\nabla w_1|^2dxdt\leq C(\overline{\n})\int |w_{10}|^2dx.\ee

Then, multiplying (\ref{sas2}) by
$ w_{1t} $ and integrating the resulting equality over $\O $  yield  that  for $t\in (0,\si(T)],$
 \be\la{sas5}\ba& \frac{1}{2}\left( \mu\|\na {w_1}\|_{L^2}^2
 +(\mu+\lambda)\|\div {w_1}\|^2_{L^2} \right)_t
 +\int\n|\dot {w_1}|^2dx\\&=\int \n\dot {w_1}(u\cdot\na {w_1})dx
   \\&\le  C(\on )\|\n^{1/2}\dot w_1\|_{L^2} \|\n^{1/(2+\nu)}u\|_{L^{2+\nu}}\| \na^2  w_1\|_{L^2}^{2/(2+\nu)} \| \na  w_1\|_{L^2}^{\nu/(2+\nu)}
  \\&\le \frac12\int\n|\dot {w_1}|^2dx+C(\on,M)  \| \na  w_1\|_{L^2}^2,
  \ea\ee
  where in the last inequality we have used \eqref{jan2} and the following simple fact:
$$ \|\na^2w_1\|_{L^2}\le C\|\n\dot w_1\|_{L^2}, $$
  due to     the standard $L^2$-estimate of the  elliptic  system  (\ref{sas2}).
 Gronwall's inequality  together with (\ref{sas5}) and (\ref{sas4})  gives
\be \label{uu1} \sup_{0\le t\le \si(T)}\|\nabla
w_1\|_{L^2}^2+\int_0^{\si(T)}\int\n|\dot {w_1}|^2dxdt\leq C(\on,M)\|\na
w_{10}\|_{L^2}^2,\ee and \be  \label{uu2} \sup_{0\le t\le
\si(T)}t\|\nabla w_1\|_{L^2}^2+\int_0^{\si(T)}t\int\n|\dot
{w_1}|^2dxdt\leq C(\on,M)\| w_{10}\|_{L^2}^2.\ee

Since the solution operator $w_{10}\mapsto w_1(\cdot,t)$ is linear,
by the   standard Stein-Weiss interpolation argument (\cite{bl}),
one can deduce from (\ref{uu1}) and (\ref{uu2}) that for any
$\theta\in [\beta,1],$\be \label{uu4} \sup_{0\le t\le
\si(T)}t^{1-\theta}\|\nabla
w_1\|_{L^2}^2+\int_0^{\si(T)}t^{1-\theta}\int\n|\dot
{w_1}|^2dxdt\leq C(\on,M)\| w_{10}\|_{\dot H^\theta}^2,\ee with a uniform
constant $C$ independent of $\theta.$

Finally, we estimate $w_2.$ It follows from a similar way as for the proof of
(\ref{h19}) and (\ref{h18}) that  \be\la{sas6}
 \|\na((2\mu+\lambda)\div w_2-P)\|_{L^2}+\|\na (\na^\perp\cdot w_2)\|_{L^2}\le C\|\n  \dot
w_2 \|_{L^2} ,\ee and that for $p\ge 2,$\be \la{sax6}\ba \|\na w_2\|_{L^p}&  \le C (\|(2\mu+\lambda)\div w_2-P\|_{L^p}+C\|P\|_{L^p}+\|\na^\perp\cdot w_2\|_{L^p}   \\&\le \de\|\n  \dot
w_2 \|_{L^2}  +C(\on,p,\de)\|\na w_2\|_{L^2}  +C(\on,p,\de)C_0^{1/p} . \ea\ee  Multiplying (\ref{sas3})
by $w_{2t}  $ and  integrating the resulting equation over $\O  $   yield that for $t\in (0,\si(T)],$
 \be \la{xac6}\ba& \frac{1}{2}\left( \mu\|\na {w_2}\|_{L^2}^2
 +(\mu+\lambda)\|\div {w_2}\|^2_{L^2} -2\int P\div w_2dx \right)_t
 +\int\n|\dot {w_2}|^2dx\\&=\int \n\dot {w_2}(u\cdot\na
 {w_2})dx-\int P_t\div w_2dx \\&\le  C(\on)\|\n^{1/2}\dot w_2\|_{L^2}\|\n^{1/(2+\nu)}u\|_{L^{2+\nu}}\| \na  w_2\|_{L^{2(2+\nu)/\nu}} -\int P_t\div w_2dx \\&\le C(\on,M)\de\|\n^{1/2}\dot w_2\|_{L^2}^2 +C(\de,\on,M)\left(\|\na w_2\|_{L^2}^2 +\|\na u\|_{L^2}^2+1\right),
  \ea\ee where in the last inequality we have used (\ref{sax6}), \eqref{jan2},  and the following simple fact:
   \bnn\ba   -\int P_t\div w_2dx     &= -\frac{1}{2\mu+\lm}\int Pu\cdot\na ((2\mu+\lm)\div w_2-P)dx \\&\quad+\frac{1}{2(2\mu+\lm)}\int P^2\div u dx \\&\le C\|Pu\|_{L^2}\|\n\dot w_2\|_{L^2}+ C\|P^2\|_{L^2}\|\na u\|_{L^2}  \\&\le \de\|\n^{1/2}\dot w_2\|_{L^2}^2 +C(\de,\on )\left( \|\na u\|_{L^2}^2+1\right),
  \ea\enn due to \eqref{a95} and \eqref{sas6}.
 Gronwall's inequality together with (\ref{xac6})  gives
\be  \label{uu3} \sup_{0\le t\le \si(T)}\|\nabla
w_2\|_{L^2}^2+\int_0^{\si(T)}\int\n|\dot {w_2}|^2dxdt\leq
C(\on,M).\ee   Taking $w_{10}=u_0 $ so that $w_1+w_2=u,$
we then derive (\ref{uv1}) from (\ref{uu4}) and (\ref{uu3}) directly.
  Thus, we finish the proof of Lemma \ref{zc1}.

We now proceed to derive a uniform (in time) upper bound for the
density, which turns out to be the key to obtain all the higher
order estimates and thus to extend the classical solution globally.
We will use an approach motivated by our previous study on the
two-dimensional Stokes approximation equations (\cite{lx}), see also \cite{hlx1}.

\begin{lemma}\la{le7}
There exists a positive constant
   $\ve_0=\ve_0 (\on ,M) $
    depending    on  $\mu ,  \lambda ,   \ga ,  a ,  \on, \beta, $ $ N_0,$ and $M$  such that,
    if  $(\n,u)$ is a smooth solution  of
   (\ref{a1})-(\ref{h2})     on $\O \times (0,T] $
   satisfying (\ref{z1}) and the assumptions in Theorem \ref{th1}, then
      \be\la{lv102}\sup_{0\le t\le T}\|\n(t)\|_{L^\infty}  \le
\frac{7\bar \n }{4}  ,\ee
      provided $C_0\le \ve_0 . $

   \end{lemma}

{\it Proof.}
  First, we rewrite the equation of the mass conservation
$(\ref{a1})_1$ as \be \la{z.3} D_t \n=g(\n)+b'(t), \ee where \bnn
D_t\n\triangleq\n_t+u \cdot\nabla \n ,\quad
g(\n)\triangleq-\frac{\n^{\ga+1 }}{2\mu+\lambda}  ,
\quad b(t)\triangleq-\frac{1}{2\mu+\lambda} \int_0^t\n Fdt. \enn

 Next, it follows from \eqref{h19}, \eqref{o3.7},   \eqref{uq2}, and \eqref{z.1} that for $t>0$ and $p\in [2,\infty),$ \be\la{z.2}\ba \|\na F(\cdot,t)\|_{L^p}&\le C(p)\|\n \dot u (\cdot,t)\|_{L^p} \\&\le C(p,\on,M)(1+t)^5\left(\|\n^{1/2}\dot u(\cdot,t)\|_{L^2}+\|\na\dot u(\cdot,t)\|_{L^2}\right),\ea\ee which, together with   the Gagliardo-Nirenberg inequality  \eqref{g2} for $q=2$, yields  that
for $r\triangleq 4+4/\beta$ and
$\de_0\triangleq (2r+(1-\beta)(r-2))/(3r-4)\in(0,1),$
\bnn\ba  &|b(\si(T)) |  \\  &\le C(\on)
\int_0^{\si(T)} \si^{(\beta-1)(r-2)/(4(r-1))}\left(\si^{1-\beta}\|F \|_{L^2}^2\right)
^{(r-2)/(4(r-1))}\|\na
F \|^{r/(2(r-1))}_{L^r}dt\\ &\le  C(\on,M) \int_0^{\si(T)}
\si^{-(2r+(1-\beta)(r-2))/(4(r-1))}\left(\si^2\|\na F\|^{2}_{L^r} \right)^{r/(4(r-1))}dt\\ &\le  C(\on,M) \left(\int_0^{\si(T)}
\si^{-\de_0}dt\right)^{(3r-4)/(4(r-1))}  \left(\int_0^{\si^(T)} \si^2\|\na F\|^{2}_{L^r}  dt\right)^{r/(4(r-1))}\\ &\le  C(\on,M) \left( \int_0^{\si(T)} \left(\si^2\|\n^{1/2}\dot u\|^{2}_{L^2} +\si^2\|\na\dot u\|^{2}_{L^2}\right)dt \right)^{r/(4(r-1))}\\ &\le  C(\on,M) C_0^{r/(4(r-1))} ,\ea\enn
where in the second, fourth, and last inequalities one has used respectively \eqref{uv1}, \eqref{z.2}, and   \eqref{h27}. This combined with \eqref{z.3} yields that
   \be\la{a103}\sup_{t\in
[0,\si(T)]}\|\n\|_{L^\infty} \le \on
+C(\on,M)C_0^{1/4}\le\frac{3 \bar\n  }{2},\ee
 provided $$C_0\le \ve_1\triangleq\min\{1, (\on/(2C(\on,M)))^{4}\}. $$

Next,  it follows from  \eqref{h19}     and \eqref{ly8}     that for $t\in[\si(T),T],$ \be \la{hg2} \ba \|F(\cdot,t)\|_{H^1} & \le C \left(\|\na u(\cdot,t)\|_{L^2}+\|P(\cdot,t)\|_{L^2}+\|\n  \dot u(\cdot,t)\|_{L^2}\right)\\ & \le C( \on) C_0^{1/2}  t^{-1/2}   , \ea\ee which together with \eqref{g2} and  \eqref{z.2} shows
\be\la{hg1}\ba &\int_{\si(T)}^T\|F(\cdot,t)\|_{L^\infty}^4dt \\&\le C\int_{\si(T)}^T\|F(\cdot,t)\|_{L^{72}}^{35/9}\|\na F(\cdot,t)\|_{L^{72}}^{1/9} dt\\ &\le C(\on,M)C_0^{35/18}\int_{\si(T)}^Tt^{-25/18}(\|\n^{1/2}\dot u\|_{L^2} +\|\na\dot u\|_{L^2} )^{1/9}dt\\ &\le C(\on,M)C_0^{35/18},\ea\ee
where in the last inequality, one has used \eqref{z1}.
This shows that for all $\si(T)\le t_1\le t_2\le T,$ \bnn\ba  |b(t_2)-b(t_1)|    &\le  C(\on) \int_{t_1}^{t_2} \|F(\cdot,t)\|_{L^\infty}dt  \\&\le \frac{1}{2\mu+\lambda}(t_2-t_1)+ C(\on,M) \int_{\si(T)}^T\|F(\cdot,t)\|_{L^\infty}^4dt\\ &\le \frac{1}{2\mu+\lambda}(t_2-t_1)+  C(\on,M) C_0^{35/18}   ,  \ea\enn which implies that  one can
choose $N_1$ and $N_0$ in (\ref{a100}) as: \bnn
N_1=\frac{ 1}{2\mu+\lambda},\quad  N_0=C(\on,M)C_0^{35/18} .\enn  Hence, we set $\bar\zeta= 1 $ in (\ref{a101}) since for all $  \zeta \ge
 1,$
$$ g(\zeta)=-\frac{ \zeta^{\ga+1}}{2\mu+\lambda} \le -N_1
=- \frac{ 1}{2\mu+\lambda}.  $$   Lemma
\ref{le1} and (\ref{a103}) thus lead to \be\la{a102} \sup_{t\in
[\si(T),T]}\|\n\|_{L^\infty}\le  \frac{ 3\bar \n
}{2}  +N_0 \le
\frac{7\bar \n }{4} ,\ee provided $$ C_0\le
\ve_0\triangleq\min\{\ve_1,\ve_2 \}, \quad\mbox{ for
}\ve_2\triangleq \left(\frac{ \bar \n }{4C(\on,M) }\right)^{18/35}.$$
The combination of (\ref{a103}) with (\ref{a102}) completes the
proof of Lemma \ref{le7}.

With Lemmas \ref{le5} and \ref{le7} at hand, we are now in a position to prove  Proposition \ref{pr1}.

{\it Proof of Proposition \ref{pr1}.}  It follows from \eqref{h27} that \be\la{lv101}  A_1(T)+A_2(T)+\int_0^T\si\|P\|_{L^2}^2dt\le C_0^{1/2},\ee provided $$C_0\le \ve_3\triangleq (C(\on))^{-2}.$$ Letting $\ve\triangleq \min\{\ve_0,\ve_3\},$ we obtain \eqref{z2} directly from \eqref{lv102} and \eqref{lv101}  and finish the proof of Proposition \ref{pr1}.

\section{\la{se5} A priori estimates (II): higher order estimates }

Form now on, for smooth initial data $(\n_0,u_0)$  satisfying \eqref{co1} and \eqref{h7},   assume that
  $(\n,u)$ is a smooth solution of (\ref{a1})-(\ref{h2}) on $\O \times (0,T] $ satisfying (\ref{z1}).
Then, we  derive some necessary uniform  estimates on the spatial gradient of
the smooth solution $(\n,u)$.

\begin{lemma}\la{le4}     There is a positive constant $C $ depending only on $T,\mu ,  \lambda ,   \ga ,  a ,  \on, \beta, N_0, M,q,$   and $\|\n_0\|_{H^1\cap W^{1,q}} $ such that
  \be\la{pa1}\ba
  &\sup_{0\le t\le T}\left(\norm[H^1\cap W^{1,q}]{ \rho} +\|\nabla   u\|_{L^2} +  t\|\na^2 u\|^2_{L^2}  \right)\\&+\int_0^T \left(\|\nabla^2  u\|_{L^2}^2+\|\nabla^2 u\|_{L^q}^{(q+1)/q}+t\|\nabla^2 u\|_{L^q}^2 \right)dt\le C .\ea
  \ee
\end{lemma}

{\it Proof.}  First, it follows from \eqref{h29}, \eqref{ly12},  Gronwall's inequality, and \eqref{a16} that
     \be\la{a93}
  \sup_{t\in[0,T]}\|\nabla u\|_{L^2}^2 + \int_0^{T}\int\n|\dot{u}|^2dxdt
  \le C,
  \ee which together with \eqref{h18} shows \be\la{hc1} \int_0^T\|\na u\|_{L^4}^4dt\le C.\ee
Multiplying \eqref{mm4} by $t$ and integrating the resulting inequality over $(0,T)$ combined with  \eqref{a93}  and \eqref{hc1}  lead to \be\la{b19}
\sup_{0\le t\le T}t\int \rho|\dot{u}|^2dx + \int_0^Tt\|\na\dot{u}\|_{L^2}^2dt\le C.
\ee

Next, we prove (\ref{pa1}) by using Lemma \ref{le9} as in
\cite{hlx}. For $   p\in [2,q],$ $|\nabla\n|^p$ satisfies \bnnn \ba
& (|\nabla\n|^p)_t + \text{div}(|\nabla\n|^pu)+ (p-1)|\nabla\n|^p\text{div}u  \\
 &+ p|\nabla\n|^{p-2}(\nabla\n)^t \nabla u (\nabla\n) +
p\n|\nabla\n|^{p-2}\nabla\n\cdot\nabla\text{div}u = 0.\ea
\ennn Thus, \be\la{L11}\ba
\frac{d}{dt} \norm[L^p]{\nabla\n}  &\le
 C(1+\norm[L^{\infty}]{\nabla u} )
\norm[L^p]{\nabla\n} +C\|\na^2u\|_{L^p}\\ &\le
 C(1+\norm[L^{\infty}]{\nabla u} )
\norm[L^p]{\nabla\n} +C\|\n\dot u\|_{L^p}, \ea\ee due to\be
\la{ua1}\|\na^2 u\|_{L^p}\le   C\left(\|\n\dot u\|_{L^p}+ \|\nabla
P \|_{L^p}\right),\ee which follows from the standard
$L^p$-estimate for the following elliptic system:
 \bnn\la{zp202}  \mu\Delta
u+(\mu+\lambda)\na {\rm div}u=\n \dot u+\na P,\quad \, u\rightarrow
0\,\,\mbox{ as } |x|\rightarrow \infty. \enn

Next, it follows from
  the  Gargliardo-Nirenberg inequality, \eqref{a93},  and  \eqref{h19}   that
 \be\la{419}\ba  \|\div u\|_{L^\infty}+\|\o\|_{L^\infty}  &\le C \|F\|_{L^\infty}+C\|P\|_{L^\infty} +C\|\o\|_{L^\infty}\\ &\le C(q) +C(q) \|\na F\|_{L^q}^{q/(2(q-1))} +C(q) \|\na \o\|_{L^q}^{q/(2(q-1))}\\ &\le C(q) +C(q) \|\n\dot u\|_{L^q}^{q/(2(q-1))} , \ea\ee
 which, together with
Lemma \ref{le9}, yields that
    \be\la{b24}\ba   \|\na
u\|_{L^\infty }  &\le C\left(\|{\rm div}u\|_{L^\infty }+
\|\o\|_{L^\infty } \right)\log(e+\|\na^2 u\|_{L^q}) +C\|\na
u\|_{L^2} +C \\&\le C\left(1+\|\n\dot u\|_{L^q}^{q/(2(q-1))}\right)\log(e+\|\rho \dot u\|_{L^q} +\|\na \rho\|_{L^q}) +C\\&\le C\left(1+\|\n\dot u\|_{L^q} \right)\log(e+   \|\na \rho\|_{L^q}) . \ea\ee

Next, it follows from the H\"older inequality   and \eqref{z.2}   that
 \bnn\la{b22}\ba
 \| \rho \dot u\|_{L^q} & \le
  \| \rho \dot u\|_{L^2}^{2(q-1)/(q^2-2)}\|\n\dot u\|_{L^{q^2}}^{q(q-2)/(q^2-2)}\\ & \le
 C \| \rho \dot u\|_{L^2}^{2(q-1)/(q^2-2)}\left(\| \rho^{1/2} \dot u\|_{L^2}+\|\na\dot u\|_{L^2}\right)^{q(q-2)/(q^2-2)}\\ & \le
 C \| \rho^{1/2}  \dot u\|_{L^2} +C \| \rho^{1/2} \dot u\|_{L^2}^{2(q-1)/(q^2-2)}\|\na \dot u\|_{L^2}^{q(q-2)/(q^2-2)}   , \ea\enn which combined with   \eqref{a93} and  \eqref{b19} implies that\be\la{4a2}   \ba &\int_0^T\left(\| \rho \dot u\|_{L^q}^{1+1 /q}+t\| \rho \dot u\|_{L^q}^2\right)  dt\\ & \le C   \int_0^T\left( \| \rho^{1/2}  \dot u\|_{L^2}^2+ t\|\na \dot u\|_{L^2}^2 + t^{-(q^3-q^2-2q-1)/(q^3-q^2-2q )}\right)dt \\ &\le  C   .\ea\ee

Then,
 substituting \eqref{b24} into \eqref{L11} where $p=q$, we deduce from Gronwall's inequality and \eqref{4a2}  that \bnn \la{b30} \sup\limits_{0\le t\le T}\|\nabla
\rho\|_{L^q}\le C,\enn  which, along with \eqref{ua1} and \eqref{4a2}, shows
  \be \la{b31}\int_0^T \left(\|\nabla^2 u\|_{L^q}^{(q+1)/q}+t\|\nabla^2 u\|_{L^q}^2\right)dt\le C.\ee

Finally,
taking $p=2$ in (\ref{L11}), one gets by using (\ref{b31}),
(\ref{a93}), and Gronwall's inequality that \bnn \sup\limits_{0\le
t\le T}\|\nabla \n\|_{L^2}\le C,\enn which, together with
(\ref{ua1}),   (\ref{b19}),  and
(\ref{b31}), yields (\ref{pa1}). The proof of Lemma \ref{le4} is
completed.

\begin{lemma} \la{le6} There is a positive  constant $C $ depending only on  $T,\mu ,  \lambda ,   \ga ,  a ,  \on, \beta, N_0,M,q,$    and $\|\na(\bar x^a\n_0)\|_{L^2\cap L^q}  $ such that
  \be \la{q} \ba
  &\sup_{0\le t\le T}   \|  \bar x^a\n \|_{L^1\cap H^1\cap W^{1,q}}  \le C .\ea
  \ee
\end{lemma}

 {\it Proof.}
First, it follows from \eqref{3h},   \eqref{o3.7}, \eqref{uq2}, and \eqref{3.15}   that for any $\eta\in(0,1]$ and any $s>2,$
 \be\la{a4.21} \|u\bar x^{-\eta}\|_{L^{s/\eta}}\le C(\eta,s).\ee
Multiplying $\eqref{a1}_1$ by $\bar x^a $ and integrating the resulting equality over $\O $ lead to
\bnn \ba  \frac{d}{dt}\int\n \bar x^a  dx &\le C\int \n |u|\bar x^{a -1} \log^2 (e+|x|^2) dx\\ &\le C \|\n\bar x^{a-1+8/(8+a)}\|_{L^{(8+a)/(7+a)}}\|u\bar x^{-4/(8+a)}\|_{L^{8+a}}\\&\le C\int\n \bar x^a  dx +C.\ea\enn
 This  gives
\be \la{o4.7} \sup_{0\le t\le T}\int\n \bar x^a  dx\le  C .\ee

Then, one derives from $\eqref{a1}_1$ that $ v\triangleq\n\bar x^a$ satisfies \bnn\ba v_t+u\cdot\na v-a vu\cdot\na \log \bar x+v\div u=0,\ea\enn
which, together with some estimates as for \eqref{L11}, gives that for any $p\in [2,q]$
\be\la{7.z1}\ba (\|\na v\|_{L^p} )_t& \le C(1+\|\na u\|_{L^\infty}+\|u\cdot \na \log \bar x\|_{L^\infty}) \|\na v\|_{L^p} \\&\quad +C\|v\|_{L^\infty}\left( \||\na u||\na\log \bar x|\|_{L^p}+\||  u||\na^2\log \bar x|\|_{L^p}+\| \na^2 u \|_{L^p}\right)\\& \le C(1 +\|\na u\|_{W^{1,q}})  \|\na v\|_{L^p} \\&\quad+C\|v\|_{L^\infty}\left(\|\na u\|_{L^p}+\|u\bar x^{-2/5}\|_{L^{4p}}\|\bar x^{-3/2}\|_{L^{4p/3}}+\|\na^2 u\|_{L^p}\right) \\& \le C(1 +\|\na^2u\|_{L^p}+\|\na u\|_{W^{1,q}})(1+ \|\na v\|_{L^p}+\|\na v\|_{L^q}), \ea\ee
where in the second and the last inequalities, one has used \eqref{a4.21} and \eqref{o4.7}.
Choosing $p=q$ in \eqref{7.z1}, we obtain after using  Gronwall's inequality and \eqref{pa1} that
\be\la{7.z2}\ba  \sup\limits_{0\le t\le T}\|\na (\n \bar x^a)\|_{L^q} \le C. \ea\ee

Finally, setting $p=2$ in  \eqref{7.z1}, we deduce from   \eqref{pa1} and   \eqref{7.z2} that
  \bnn \sup\limits_{0\le t\le T}\|\na(\n \bar x^a)\|_{L^2 } \le C , \enn
 which combined with \eqref{o4.7} and \eqref{7.z2} thus gives
  \eqref{q} and finishes
 the proof of Lemma \ref{le6}.

\begin{lemma}\label{lem4.5v}  There is a positive  constant $C $ depending only on  $T,\mu ,  \lambda ,   \ga ,  a ,  \on, \beta, N_0,M,q,$    and $\|\na(\bar x^a\n_0)\|_{L^2\cap L^q}  $ such that
\begin{equation}\la{5.13v}\ba
 \sup_{0\leq t\leq T }t\left(\|\n^{1/2}u_t\|^2_{L^2}+\|\na u\|^2_{H^1} \right)+\int_0^Tt\|\na u_t\|_{L^2}^2dt\le C.\ea
\end{equation}
\end{lemma}

{\it Proof.}
Differentiating $\eqref{a1}_2$ with respect to $t$ gives
\be\la{zb1}\ba &\n u_{tt}+\n u\cdot \na u_t-\mu\Delta u_t-( \mu+\lm)\na  \div u_t  \\ &=-\n_t(u_t+u\cdot\na u)-\n u_t\cdot\na u -\na P_t.\ea\ee
 Multiplying \eqref{zb1} by $u_t$ and integrating the resulting equation over $\O,$ we obtain after using  $\eqref{a1}_1$ that\be\ba  \la{na8}&\frac{1}{2}\frac{d}{dt} \int \n |u_t|^2dx+\int \left(\mu|\na u_t|^2+( \mu+\lm)(\div u_t)^2  \right)dx\\
 &=-2\int \n u \cdot \na  u_t\cdot u_tdx  -\int \n u \cdot\na (u\cdot\na u\cdot u_t)dx\\
  &\quad-\int \n u_t \cdot\na u \cdot  u_tdx
 +\int P_{t}\div u_{t} dx\\
 &\le C\int  \n |u||u_{t}| \left(|\na  u_t|+|\na u|^{2}+|u||\na^{2}u|\right)dx +C\int \n |u|^{2}|\na u ||\na u_{t}|dx \\
 &\quad+C\int \n |u_t|^{2}|\na u |dx +C(\de)\|P_{t}\|_{L^2}^2+\de\|\na u_{t}\|_{L^2}^2 .
  \ea\ee

  Each term  on the right-hand side of  \eqref{na8} can be estimated as follows:

First, the combination of \eqref{a4.21} with \eqref{q} gives  that
for any $\eta\in(0,1]$ and any $s>2,$
 \be\la{5.d2}\|\n^\eta u \|_{L^{s/\eta}}+ \|u\bar x^{-\eta}\|_{L^{s/\eta}}\le C(\eta,s).\ee Moreover,
it follows from    \eqref {z.1}, \eqref{o3.7}, and \eqref{uq2}   that   \be \la{5.d1}\|\n^{1/2} u_t\|_{L^6} \le C \|\n^{1/2}u_t\|_{L^2}+C \|\na u_t\|_{L^2},\ee
 which together with \eqref{5.d2}, \eqref{a93}, and Holder's inequality  yields  that for $\de\in (0,1),$
 \be\la{na2}\ba  &\int  \n |u||u_{t}| \left(|\na  u_t|+|\na u|^{2}+|u||\na^{2}u|\right) dx\\
 & \le C \|\n^{1/2} u\|_{L^{6}}\|\n^{1/2} u_{t}\|_{L^{2}}^{1/2} \|\n^{1/2} u_{t}\|_{L^{6}}^{1/2}\left(\| \na u_{t}\|_{L^{2}}+\| \na u\|_{L^{4}}^{2} \right) \\
 &\quad +C\|\n^{1/4}  u \|_{L^{12}}^{2}\|\n^{1/2} u_{t}\|_{L^{2}}^{1/2} \|\n^{1/2} u_{t}\|_{L^{6}}^{1/2} \| \na^{2} u \|_{L^{2}}  \\
 & \le C  \|\n^{1/2} u_{t}\|_{L^{2}}^{1/2}\left(\|\n^{1/2} u_{t}\|_{L^{2}} +\| \na u_{t}\|_{L^{2}}\right)^{1/2}\left(\| \na u_{t}\|_{L^{2}}+  \| \na^{2} u \|_{L^{2}}+1\right)\\
 &\le  \de\| \na u_{t}\|_{L^{2}}^{2}+C(\de)  \left(\| \na^{2} u \|_{L^{2}}^{2} +  \|\n^{1/2} u_{t}\|_{L^{2}}^{2}+1\right).
 \ea\ee

 Next, Holder's inequality, \eqref{5.d2},  and \eqref{5.d1} lead to
 \be \ba \la{5.ap3}&  \int \n |u|^{2}|\na u ||\na u_{t} |dx+\int \n |u_t|^{2}|\na u |dx \\ &\le C \|\n^{1/2} u\|_{L^{8}}^{2}\|\na u\|_{L^{4}} \| \na u_{t}\|_{L^{2}}+\| \na u\|_{L^{2}}
   \|\n^{1/2} u_{t}\|_{L^{6}}^{3/2}\|\n^{1/2} u_{t}\|_{L^{2}}^{1/2} \\
 &\le  \de\| \na u_{t}\|_{L^{2}}^{2}+C(\de)  \left(\| \na^{2} u \|_{L^{2}}^{2} +  \|\n^{1/2} u_{t}\|_{L^{2}}^{2}+1\right).\ea\ee

Next, it follows from \eqref{5.d2},  \eqref{a93},  and \eqref{q} that
\be \la{nan3}\ba \|P_t\|_{L^2 } &\le C\|\bar x^{-a} u\|_{L^{2q/(q-2)}}\|\n\|_{L^\infty}^{\ga-1}\|\bar x^a \na \n\|_{  L^q}+C\|\na u\|_{L^2 }\le C.\ea\ee

Finally, putting
  \eqref{na2}-\eqref{nan3} into \eqref{na8} and choosing $\de$ suitably  small, we obtain after  using \eqref{ua1} and \eqref{pa1} that
\be \la{5.av2} \frac{d}{dt} \int \n |u_t|^2dx+\mu\int   |\na u_t|^2 dx   \le   C \int \n |u_t|^2dx+C\int \n |\dot u |^2dx+1.
\ee  It follows from \eqref{ua1}  and \eqref{5.d2} that
  \bnn\ba  & \|\na u\|_{H^1}+\|\n^{1/2} u\cdot\na u\|_{L^2}\\&\le C+C\|\n^{1/2}  \dot u\|_{L^2}+ C\|\n^{1/2} u\|_{L^6}\|\na u\|_{L^2}^{2/3}\|\na^2u\|_{L^2}^{1/3}\\&\le C +C\|\n^{1/2}  \dot u\|_{L^2}+\frac12\|\na^2u\|_{L^2},\ea\enn which together with  \eqref{a93}   shows
  \be \la{5.av3} \|\na u\|_{H^1}+\|\n^{1/2} u_t\|_{L^2}\le C\|\n^{1/2} \dot u\|_{L^2}+C.\ee
  This combined   with   \eqref{5.av2}, \eqref{a93},    and Gronwall's inequality gives \eqref{5.13v} and finishes the proof of Lemma \ref{lem4.5v}.

From now on,  assume that
  $(\n,u)$ is a smooth solution of (\ref{a1})-(\ref{h2}) on $\O \times (0,T] $ satisfying (\ref{z1}) for smooth initial data $(\n_0,u_0)$  satisfying \eqref{co1},  \eqref{h7},  \eqref{1.c1},  and \eqref{co2}.
Moreover,   in addition to    $T, \mu ,  \lambda ,   \ga ,  a ,  \on, \beta,  N_0,M,q,$   and  $\|\na(\bar x^a\n_0)\|_{L^2\cap L^q}  ,$ the generic positive
constant $C $ may depend on   $\|\na^2u_0\|_{L^2},$   $\| \bar x^{\de_0} \na^2 \rho_0 \|_{L^2},$   $\|\bar x^{\de_0}  \na^2  P(\rho_0) \|_{L^2},\de_0,$  and $\|g\|_{L^2},$
 with $g$   as in (\ref{co2}).

\begin{lemma}\label{lem4.5}  It holds that
\begin{equation}\la{5.13g}\ba
 \sup_{0\leq t\leq T }\left(\|\n^{1/2}u_t\|_{L^2}+\|\na u\|_{H^1} \right)+\int_0^T\|\na u_t\|_{L^2}^2dt\le C.\ea
\end{equation}
\end{lemma}

{\it Proof.} Taking into account on the compatibility condition (\ref{co2}), we can
  define \be \la{5.av1}  \sqrt{\n} \dot u (x,t=0)= g  .\ee
  Integrating \eqref{mm4}  over $(0,T)$ together with \eqref{5.av1} and   \eqref{hc1}   yields directly that
  \bnn \la{5.11v}\ba
 \sup_{0\leq t\leq T } \|\n^{1/2}\dot u \|_{L^2} +\int_0^T\|\na\dot u \|_{L^2}^2dt\le C,\ea\enn
 which, along with \eqref{5.av3} and  \eqref{5.av2},    gives
   \eqref{5.13g} and finishes the proof of Lemma \ref{lem4.5}.

   The following
higher order estimates of the solutions which are needed to
guarantee the extension of local classical solution to be a global
one  are similar to those in \cite{hlma}, so we omit  their  proofs here.
 \begin{lemma} \la{4.m44} The following estimates hold:
\begin{equation} \ba
 \sup_{0\leq t\leq T}\left(\|\bar x^{ \de_0}\na^2\n\|_{L^2}+ \|\bar x^{ \de_0}\na^2 P \|_{L^2}\right)\le C,\ea
\end{equation}
\be\ba \la{oo8}\sup\limits_{0\le t\le T}t \|\na u_t\|_{L^{2}}^2 +\int_{0}^{T}t\left(\|\n^{1/2}  u_{tt}\|_{L^{2}}^2+\|\na^2  u_{t}\|_{L^{2}}^2\right)dt \le C ,\ea\ee
\begin{equation}\la{5.13n}\ba
 \sup_{0\leq t\leq T}\left(\|\nabla^2 \n\|_{L^q } +\|\nabla^2 P  \|_{L^q }\right) \leq C ,\ea
\end{equation}
\begin{equation}\label{4.m30}\ba
&\sup_{0\leq t\leq T} t\left(\|\n^{1/2}u_{tt}\|_{L^2}+  \|\na^3
u \|_{L^2\cap L^q} + \|\na
u_t \|_{H^1} +\|\na^2(\n u)\|_{L^{(q+2)/2}} \right)\\& +\int_{0}^T  t^2\left( \|\nabla u_{tt}\|_{L^2}^2 +\|u_{tt}\bar x^{-1}\|_{L^2}^2 \right)dt\leq
C .\ea
\end{equation}
 \end{lemma}

\section{\la{se4}Proofs of  Theorems  \ref{th1}-\ref{thv}}

With all the a priori estimates in Sections \ref{se3} and \ref{se5} at hand, we are
ready to prove the main results of this paper in this section.

{\it Proof of Theorem \ref{th1}.} By Lemma \ref{th0}, there exists a
$T_*>0$ such that the Cauchy problem (\ref{a1})-(\ref{h2}) has a unique strong solution $(\n,u)$ on $\O \times
(0,T_*]$. We will use the a priori estimates, Proposition \ref{pr1}
and Lemmas \ref{le4}-\ref{lem4.5v}, to extend the local strong
solution $(\n,u)$ to all time.

First,  it follows from (\ref{As1}), (\ref{As2}),  and
(\ref{co1})    that
$$ A_1(0)+A_2(0)=0,   \quad \n_0\leq
\bar{\n}.$$   Therefore, there exists a
$T_1\in(0,T_*]$ such that (\ref{z1}) holds for $T=T_1$.

Next, set \bn \la{s1}T^*=\sup\{T\,|\,{\rm (\ref{z1}) \
holds}\}.\en Then $T^*\geq T_1>0$. Hence, for any $0<\tau<T\leq T^*$
with $T$ finite, one deduces from   \eqref{5.13v}
that for any $q\ge 2,$
 \be \la{sp43}\na u \in C([\tau ,T];L^2\cap L^q) ,\ee where one has used the standard
embedding
$$L^\infty(\tau ,T;H^1)\cap H^1(\tau ,T;H^{-1})\hookrightarrow
C\left([\tau ,T];L^q\right),\quad\mbox{ for any } q\in [2,\infty).  $$
Moreover, it follows from \eqref{pa1}, \eqref{q}, and \cite[Lemma 2.3]{L2} that \be \la{n20}\n\in C([0,T];L^1\cap H^1\cap W^{1,q}).\ee

Finally,   we claim that \be \la{s2}T^*=\infty.\ee Otherwise,
$T^*<\infty$. Then by Proposition \ref{pr1}, (\ref{z2}) holds for
$T=T^*$. It follows from  \eqref{a16}, \eqref{q},  \eqref{sp43} and
(\ref{n20}) that $(\n(x,T^*),u(x,T^*))$ satisfies
(\ref{co1})   except $ u(\cdot,T^*)\in \dot H^\beta.$
  Thus, Lemma
\ref{th0} implies that there exists some $T^{**}>T^*$, such that
(\ref{z1}) holds for $T=T^{**}$, which contradicts (\ref{s1}).
Hence, (\ref{s2}) holds. Lemmas \ref{th0} and \ref{le4}-\ref{lem4.5v}
  thus show that $(\n,u)$ is in fact the unique
strong  solution defined on $\O \times(0,T]$ for any
$0<T<T^*=\infty$. The proof of Theorem  \ref{th1} is completed.

{\it Proof of Theorem \ref{t1}.} Similar to the   proof of Theorem  \ref{th1}, one can prove     Theorem  \ref{t1} by using Lemma \ref{th0}, Proposition \ref{pr1},
and Lemmas \ref{le4}-\ref{4.m44}.

To prove Theorem  \ref{thv}, we need the following  elementary estimates similar to those of Lemma \ref{le3} whose proof can be found in \cite[Lemma 2.3]{hlx1}.
\begin{lemma} \la{lve3}
  Let $\Omega=\r^3$ and $(\rho,u)$ be a smooth solution of
   (\ref{a1}).
    Then there exists a generic positive
   constant $C$ depending only on $\mu$ and $\lambda$ such that for any $p\in [2,6]$
          \be\la{hv19}
    \|{\nabla F}\|_{L^p} + \|{\nabla \o}\|_{L^p}
   \le C\norm[L^p]{\rho\dot{u}},\ee \be
      \la{hv20}\norm[L^p]{F} + \norm[L^p]{\o}
   \le C \norm[L^2]{\rho\dot{u}}^{(3p-6)/(2p) }
   \left(\norm[L^2]{\nabla u}
   + \norm[L^2]{P }\right)^{(6-p)/(2p)} ,
\ee \be \la{hv18}
   \norm[L^p]{\nabla u} \le C \left(\norm[L^p]{F} + \norm[L^p]{\o}\right)+
   C \norm[L^p]{P },
  \ee where $F=(2\mu+\lambda)\div u-P$ and $\o=\na\times u$ are the effective viscous flux and the vorticity respectively.
\end{lemma}

{\it Proof of Theorem  \ref{thv}.} It suffices to prove \eqref{lvy8}. In fact, it follows from \cite[Proposition 3.1 and (3.6)]{hlx1} that there exists some $\ve$ depending only on $\mu,\lambda,\gamma,\bar\n,\beta,$ and $M$  such that \be \la{hv27}\ba&\sup\limits_{1\le t<\infty}\left( \|\na u\|_{L^2} +\|\n\|_{L^\ga\cap L^\infty}+ \|\n^{1/2}\dot u\|_{L^2}\right)\\&+\int_1^\infty\left(\|\na u\|_{L^2}^2+ \|\n^{1/2}\dot u\|_{L^2}^2+ \|\na \dot u\|_{L^2}^2\right)dt\le C,\ea\ee provided $C_0\le \ve.$

If $\ga\le 3/2,$ it then holds that \be \la{vo1}\sup\limits_{1\le t<\infty}\|\n\|_{L^{3/2}}\le C\sup\limits_{1\le t<\infty}\|\n\|_{L^\ga}^{2\ga/3} \le C.\ee
If $\ga> 3/2,$ since $\n_0\in L^1,$ $\eqref{a1}_1$ yields that for $t\ge 0,$
\bnn \int \n(x,t)dx=\int\n_0(x)dx,\enn which combined with
\eqref{hv27} implies \be \la{vo2}\sup\limits_{1\le t<\infty}\|\n\|_{L^{3/2}}\le C\sup\limits_{1\le t<\infty}\|\n\|_{L^1}^{2 /3} \le C.\ee

 Similar to \eqref{u}, one deduces from $(\ref{a1})_2$ that
\bnn P=(-\Delta )^{-1}\div(\n  \dot u )+(2\mu+\lambda)\div u,\enn
which together with the Sobolev inequality gives
\bnn \ba \|P\|_{L^2} \le&C \|(-\Delta)^{-1} \div(\n  \dot u )\|_{L^{2}}   +C\|\na u\|_{L^2}  \\ \le& C\|  \n  \dot u \|_{L^{6/5}}   +C\|\na u\|_{L^2} \\ \le& C\|\n\|_{L^{3/2}}^{1/2}\| \n^{1/2}  \dot u  \|_{L^2}  +C\|\na u\|_{L^2} \\ \le& C \| \n^{1/2}  \dot u  \|_{L^2} +C\|\na u\|_{L^2} ,\ea \enn where in the last inequality one has used \eqref{vo1} and  \eqref{vo2}.
 This combined with \eqref{hv27} leads to
\be \la{av16}\int_1^\infty \|P\|_{L^2}^2dt\le C. \ee

Next, similar to \eqref{vv1}, for $p\ge 2,$ we have
 \be\la{vv2}\ba
\left(\| P\|_{L^p}^p \right)_t+ \frac{p\ga-1}{2\mu+\lambda}\|P\|_{L^{p+ 1}}^{p+1} &
=- \frac{p\ga-1}{2\mu+\lambda}\int
P^pFdx ,\ea\ee
  which together with Holder's inequality yields
\be\la{av96}\ba
\left(\| P\|_{L^p}^p \right)_t+\frac{p\ga-1}{2(2\mu+\lambda)} \|P\|_{L^{p+1}}^{p+1}
    &\le C(p) \|F\|_{L^{p+1}}^{p+1}  .\ea\ee

Next, for $B(t)$ defined as in \eqref{nv1}, it follows from \eqref{n1} and \eqref{hv18} that
\be   \la{lvy11}\ba
  B' (t)  + \int  \n|\dot{u} |^2dx & \le
   C  \| P\|_{L^3}^3
 + C  \|\nabla u\|_{L^3}^3   \\& \le
  C_1 \|P\|_{L^3}^3+ C\|F\|_{L^3}^3+ C\|\o\|_{L^3}^3.\ea\ee
Choosing $ C_2\ge 2+ {2(2\mu+\lm)(C_1+1)}/(2\ga-1)$
suitably large such that \be\la{lvy12} \frac{\mu}{4}\|\na u\|_{L^2}^2+\|P\|_{L^2}^2\le B(t)+C_2\|P\|_{L^2}^2\le C\|\na u\|_{L^2}^2+C\|P\|_{L^2}^2,\ee   setting $p=2 $ in \eqref{av96}, and  adding \eqref{av96} multiplied by $C_2$ to \eqref{lvy11}  yield that for $t\ge 1,$
\be \la{hv29}\ba &
 2\left( B (t)+  C_2 \| P\|_{L^2}^2\right)' +    2\int  \left(\n|\dot{u} |^2 + P^3  \right)dx  \\& \le  C\|F\|_{L^3}^3+ C\|\o\|_{L^3}^3\\&\le
   \|\n^{1/2}\dot u\|_{L^2}^2+C \left(\|\na u\|_{L^2}^4+\|P\|_{L^2}^4\right),\ea\ee where in the second inequality we have used \eqref{hv20} and \eqref{hv27}. Multiplying \eqref{hv29} by $t,$  along with Gronwall's inequality,   \eqref{lvy12},    \eqref{hv27}, and  \eqref{av16},   gives
 \be \la{hv31} \ba & \sup\limits_{
1\le t<\infty}t\left( \| \na u\|_{L^2}^2+ \| P\|_{L^2}^2
 \right)  + \int_{1}^\infty t \int  \left(\n|\dot{u} |^2 + P^3  \right)dxdt\le C . \ea\ee

Then,
multiplying \eqref{mm4} by $t^2 $ together with \eqref{hv18} gives \be\la{zvo1} \ba & \left(t^2
  \int\n|\dot{u}|^2dx \right)_t + { \mu}t^2 \int
 |\nabla\dot{u}|^2dx \\&
\le   2t   \int\n|\dot{u}|^2dx+ C  t^2\|F\|_{L^4}^4+ C  t^2\|\o\|_{L^4}^4 + \ti C  t^2\|P\|_{L^4}^4.\ea\ee
  Setting $p=3 $ in  \eqref{av96} and  adding \eqref{av96}   multiplied by $  2(2\mu+\lambda)(\ti C +1) t^2/( 3\ga-1 ) $ to \eqref{zvo1} lead  to
\bnn\la{lvy7} \ba & \left(t^2
  \int\n|\dot{u}|^2dx +  \frac{2(2\mu+\lambda)(\ti C +1)}{3\ga-1}t^2\| P\|_{L^3}^3\right)_t + { \mu}t^2 \|\nabla\dot{u}\|^2_{L^2}+t^2\|P\|_{L^4}^4 \\&
\le   C   t \int\left(\n|\dot{u}|^2+   P^3\right)dx+ C t^2\|F \|_{L^4}^4 + C t^2\|\o \|_{L^4}^4 \\&
\le   C   t \int\left(\n|\dot{u}|^2+   P^3\right)dx+ C t^2\|\n^{1/2} \dot u\|_{L^2}^3 \left(\|\na u
 \|_{L^2} +  \|P
 \|_{L^2} \right)  \\&
\le   C   t \int\left(\n|\dot{u}|^2+   P^3\right)dx+ C t^2\|\n^{1/2} \dot u\|_{L^2}^2 \left(\|\n^{1/2} \dot u\|_{L^2}^2 +\|\na u
 \|_{L^2}^2 +  \|P
 \|_{L^2}^2 \right),\ea\enn where in the second inequality we have used \eqref{hv20}. This combined with  Gronwall's inequality, \eqref{hv31},  \eqref{hv27}, and  \eqref{av16} yields
\be\la{lvy9}   \ba   \sup\limits_{1\le t <\infty} t^2\int  \left(\n|\dot{u} |^2 + P^3  \right)dx +  \int_{1}^\infty t^2 \left(\|\nabla\dot{u}\|^2_{L^2}+ \|P\|_{L^4}^4\right)dt \le C . \ea\ee
This combined   with \eqref{h18}   gives \eqref{lvy8} provided we show that for   $m=1,2,\cdots,
$  \be \la{lvy10}\sup\limits_{1\le t<\infty }t^{m}\|P\|_{L^{m+1}}^{m+1}+\int_0^\infty t^m\|P\|_{L^{m+2}}^{m+2}dt\le C(m) ,\ee which will be proved by induction. Since \eqref{hv31} shows that \eqref{lvy10} holds for $m=1,$ we assume that \eqref{lvy10} holds for $m=n,$ that is, \be \la{lvy16}\sup\limits_{1\le t<\infty }t^{n}\|P\|_{L^{n+1}}^{n+1}+\int_1^\infty t^n\|P\|_{L^{n+2}}^{n+2}dt\le C(n) .\ee
Setting  $p=n+2 $ in \eqref{vv2} and multiplying \eqref{vv2}  by $t^{n+1} $ give
\be\la{zvo8}\ba &
  \frac{2(2\mu+\lambda)}{(n+2)\ga-1}\left(t^{n+1}\| P\|_{L^{n+2}}^{n+2} \right)_t+ t^{n+1}\|P\|_{L^{n+3}}^{n+3} \\ &
\le C (n) t^{n}\| P\|_{L^{ n+2} }^{n+2} +C(n)   t^{n+1}\|P\|^{n+2}_{L^{n+2}} \|F\|_{L^\infty}. \ea\ee
It follows from the Gagliardo-Nirenberg  inequality, \eqref{hv19}, and \eqref{lvy9}  that
\bnn \ba\int_1^\infty  \|F\|_{L^\infty}dt
  &\le C\int_1^\infty  \|F\|_{L^6}^{1/2} \|\na F\|_{L^6}^{1/2}dt\\ &\le C\int_1^\infty   \|\n \dot u\|_{L^2}^{1/2} \|\n \dot u\|_{L^6}^{1/2} dt\\ &\le C\int_1^\infty  t^{-1/2} \|\na \dot u\|_{L^2}^{1/2}dt\\ &\le C,\ea\enn which, along   with \eqref{zvo8}, \eqref{lvy16}, and Gronwall's inequality, thus   shows that \eqref{lvy10} holds for $m=n+1.$   By induction, we obtain \eqref{lvy10} and finish the proof of  \eqref{lvy8}.  The proof of Theorem  \ref{thv} is completed.


\begin{thebibliography}{99}

\bibitem{bkm} Beale, J. T.; Kato, T.; Majda. A.
Remarks on the breakdown of smooth solutions for the 3-D Euler
equations.  \textit{Comm. Math. Phys.}   \textbf{94} (1984), 61-66.

\bibitem{bl} Bergh, J.;   Lofstrom, J.  \textit{Interpolation spaces, An
introduction.} Springer-Verlag, Berlin-Heidelberg-New York,  1976.



\bibitem{K1} Cho, Y.; Choe, H. J.;   Kim, H.
Unique solvability of the initial boundary value problems for
compressible viscous fluid. \textit{J. Math. Pures Appl.}  \textbf{83}
 (2004), 243-275.



\bibitem{K3} Cho, Y.;   Kim, H.
On classical solutions of the compressible Navier-Stokes equations
with nonnegative initial densities. \textit{Manuscript Math.}  \textbf{ 120}
 (2006), 91-129.


\bibitem{K2} Choe, H. J.;    Kim, H.
Strong solutions of the Navier-Stokes equations for isentropic
compressible fluids. \textit{J. Differ. Eqs.}  \textbf{190} (2003), 504-523.



\bibitem{F2}  Feireisl, E.  Dynamics of viscous compressible fluids. Oxford
University Press, New York,   2004.

\bibitem{F1} Feireisl, E.; Novotny, A.; Petzeltov\'{a}, H.  On the existence of globally defined weak solutions to the
Navier-Stokes equations. \textit{J. Math. Fluid Mech.}  \textbf{3} (2001),  no. 4, 358-392.

\bibitem{la} Gilbarg  D.;  Trudinger  N. S. Elliptic partial differential equations of second order. Second edition. Springer-Verlag, Berlin (1983).

\bibitem{Hof} Hoff, D.
Global existence for 1D, compressible, isentropic Navier-Stokes
equations with large initial data. \textit{Trans. Amer. Math. Soc.} \textbf{
303} (1987), no. 1, 169-181.

\bibitem{H3}Hoff, D. Global solutions of the Navier-Stokes equations
 for multidimensional compressible flow with discontinuous initial data.
\textit{J. Differ. Eqs.}   \textbf{120} (1995), no. 1, 215-254.

\bibitem{Hof2}Hoff, D.
Strong convergence to global solutions for multidimensional flows of
compressible, viscous fluids with polytropic equations of state and
discontinuous initial data.  \textit{Arch. Rational Mech. Anal.}   \textbf{132} (1995), 1-14.



\bibitem{hof2002}Hoff, D.  Dynamics of singularity surfaces for compressible,
viscous flows in two space dimensions. \textit{Comm. Pure Appl. Math.} \textbf{55}(2002), no. 11,  1365-1407.


\bibitem{Ho3}Hoff, D.  Compressible flow in a half-space with Navier boundary
  conditions. \textit{J. Math. Fluid Mech.}  \textbf{7} (2005), no. 3, 315-338.



\bibitem{hx2} Huang, X. D.; Li, J.; Xin Z. P.
Blowup criterion for viscous barotropic flows with vacuum states.
\textit{Comm. Math. Phys.}  \textbf{301} (2011), no. 1, 23-35.

\bibitem{hlx} Huang, X. D.; Li, J.; Xin Z. P.
Serrin type criterion for the three-dimensional compressible flows.
 \textit{SIAM J. Math. Anal.},   43 (2011), no. 4, 1872�C1886.


\bibitem{hlx1} Huang, X. D.; Li, J.; Xin, Z. P.    Global well-posedness of classical solutions with large
oscillations and vacuum to the three-dimensional isentropic
compressible Navier-Stokes equations.   Comm. Pure Appl.
Math.  {\bf65}, 549--585 (2012)



\bibitem{kato}Kato, T. Remarks on the Euler and Navier-Stokes equations in $R^2$. Proc. Symp. Pure Math. Vol. 45, Amer. Math. Soc., Providence, 1986, 1-7.

\bibitem{Kaz} Kazhikhov, A. V.;  Shelukhin, V. V.
Unique global solution with respect to time of initial-boundary
value problems for one-dimensional equations of a viscous gas.
\textit{Prikl. Mat. Meh.}  \textbf{41}   (1977), 282-291.


\bibitem{hlma}   Li, J.; Liang,   Z.  On    classical solutions   to the Cauchy problem of the two-dimensional barotropic
compressible Navier-Stokes equations with vacuum. http://arxiv.org/abs/1306.4752

\bibitem{lx}Li, J.; Xin, Z.  Some uniform estimates and blowup behavior of
global strong solutions to the Stokes approximation equations for
two-dimensional compressible flows. \textit{J. Differ. Eqs.}   \textbf{221} (2006), no. 2,  275-308.


\bibitem{L2} Lions,  P. L.  {Mathematical topics in fluid
mechanics. Vol. {\bf 1}. Incompressible models.}  Oxford
University Press, New York, 1996.

\bibitem{L1} Lions, P. L.    {Mathematical topics in fluid mechanics. Vol. 2. Compressible models.}  Oxford
University Press, New York,  1998.


\bi{hlx3} Z. Luo. Global existence of classical solutions to
two-dimensional Navier-Stokes equations
with Cauchy data containing vacuum, Math. Methods Appl. Sci., in press. DOI: 10.1002/mma.2896.



\bibitem{M1} Matsumura, A.;  Nishida, T.   The initial value problem for the equations of motion of viscous and heat-conductive
gases. \textit{J. Math. Kyoto Univ.}  \textbf{20}(1980), no. 1, 67-104.



\bibitem{Na} Nash, J.  Le probl\`{e}me de Cauchy pour les \'{e}quations
diff\'{e}rentielles d'un fluide g\'{e}n\'{e}ral. \textit{Bull. Soc. Math.
France.}  \textbf{90} (1962), 487-497.


\bibitem{nir}
 Nirenberg, L.  On elliptic partial differential equations.  Ann. Scuola Norm. Sup. Pisa (3), {\bf 13}, 115--162 (1959)


\bibitem{R} Rozanova, O.   Blow up of smooth solutions to the compressible
Navier-Stokes equations with the data highly decreasing at infinity.
\textit{J. Differ. Eqs.}   \textbf{245} (2008),  1762-1774.

\bibitem{S2} Salvi, R.;  Straskraba, I.
Global existence for viscous compressible fluids and their behavior
as $t\rightarrow \infty$. \textit{J. Fac. Sci. Univ. Tokyo Sect. IA. Math.}
 \textbf{40} (1993), 17-51.

\bibitem{Ser1} Serre, D.
Solutions faibles globales des \'equations de Navier-Stokes pour un
fluide compressible. \textit{C. R. Acad. Sci. Paris S\'er. I Math.}
  \textbf{303} (1986), 639-642.

\bibitem{Ser2} Serre, D.
Sur l'\'equation monodimensionnelle d'un fluide visqueux,
compressible et conducteur de chaleur. \textit{C. R. Acad. Sci. Paris S\'er.
I Math.}  \textbf{303} (1986), 703-706.

\bibitem{se1} Serrin, J.  On the uniqueness of compressible fluid motion.
\textit{Arch. Rational. Mech. Anal.}  \textbf{3} (1959), 271-288.


\bibitem{X1} Xin, Z. P.
Blowup of smooth solutions to the compressible {N}avier-{S}tokes
equation with compact density. \textit{Comm. Pure Appl. Math.}    \textbf{51}
 (1998), 229-240.

\bibitem{xy1}Xin, Z. P.; Yan, W. On blowup of classical solutions to the compressible Navier-Stokes equations. Comm. Math. Phys. \textbf{321} (2013), no. 2, 529-541.

\bibitem{zl1}Zlotnik, A. A.    Uniform estimates and stabilization of symmetric
solutions of a system of quasilinear equations.  \textit{Diff. Eqs.}
  \textbf{36} (2000),  701-716.
 \end{thebibliography}
\end{document}